\documentclass[preprint,12pt]{elsarticle}

\usepackage{amssymb,latexsym,amsmath} %Standard packages 
\usepackage{graphicx} 
\usepackage{caption}
\usepackage{subcaption}
\usepackage{textcomp}
\usepackage{color}
\journal{Journal of Computational Physics}

\begin{document}

\begin{frontmatter}

%% Title, authors and addresses

%% use the tnoteref command within \title for footnotes;
%% use the tnotetext command for theassociated footnote;
%% use the fnref command within \author or \address for footnotes;
%% use the fntext command for theassociated footnote;
%% use the corref command within \author for corresponding author footnotes;
%% use the cortext command for theassociated footnote;
%% use the ead command for the email address,
%% and the form \ead[url] for the home page:
%% \title{Title\tnoteref{label1}}
%% \tnotetext[label1]{}
%% \author{Name\corref{cor1}\fnref{label2}}
%% \ead{email address}
%% \ead[url]{home page}
%% \fntext[label2]{}
%% \cortext[cor1]{}
%% \address{Address\fnref{label3}}
%% \fntext[label3]{}

\title{Basis adaptation and domain decomposition for steady partial differential equations with random coefficients}

%% use optional labels to link authors explicitly to addresses:
%% \author[label1,label2]{}
%% \address[label1]{}
%% \address[label2]{}

\author{Ramakrishna Tipireddy, Panos Stinis and Alexandre Tartakovsky}

\address{Pacific Northwest National Laboratory, Richland WA 99354}

\begin{abstract}
We present a novel approach for solving steady-state stochastic partial differential equations (PDEs) with high-dimensional random parameter space. The proposed approach combines spatial domain decomposition with basis adaptation for each subdomain. The basis adaptation is used to address the curse of dimensionality by constructing an accurate  low-dimensional representation of the stochastic PDE solution (probability density function and/or its leading statistical moments) in each subdomain. Restricting the basis adaptation to a specific subdomain affords finding a locally accurate solution. Then, the solutions from all of the subdomains are stitched together to provide a global solution. We support our construction with numerical experiments for a steady-state diffusion equation with a random spatially dependent coefficient. Our results show that highly accurate global solutions can be obtained with significantly reduced computational costs.
\end{abstract}

\begin{keyword}
basis adaptation \sep dimension reduction \sep domain decomposition \sep polynomial chaos \sep uncertainty quantification

%% keywords here, in the form: keyword \sep keyword

%% PACS codes here, in the form: \PACS code \sep code

%% MSC codes here, in the form: \MSC code \sep code
%% or \MSC[2008] code \sep code (2000 is the default)

\end{keyword}

\end{frontmatter}

%% \linenumbers

%% main text
\section{Introduction}
\label{RK:sec:intro}

Uncertainty quantification for systems with a large number of deterministically unknown input parameters is a formidable computational task. In the probabilistic framework, uncertain parameters are treated as random variables/fields, yielding the governing equations stochastic. The two most popular methods for solving stochastic equations are Monte Carlo (MC) and polynomial chaos (PC) expansion. In both methods, the random input parameters are represented with $d$ random variables using truncated Karhunen-Lo\`eve (KL) expansion  \cite{RK:Loeve1977}. MC methods are robust and easy to implement, but they converge at a very slow rate. Hence, they require a large number of samples. On the other hand, the MC convergence rate does not depend on $d$. Contrary to this, the computational cost of standard PC methods increases exponentially with increasing $d$, a phenomenon often referred to as the ``curse of dimensionality'' \cite{RK:Ghanem1991, RK:Xiu2002, Babuka2002, Babuka2010}. Because of this, standard PC methods are only efficient for small to moderate $d$ \cite{Lin2009AWR,Lin2010JSC,Venturi2013JCP}.

Recently, basis adaptation \cite{RK:Tipireddy2014} and active subspace methods \cite{RK:Constantine2014} have been presented to identify a low-dimensional representation of the solution of stochastic equations. In the current work, we propose a new approach that combines basis adaptation with spatial decomposition (\cite{RK:Toselli2005, Agarwal20097662,Lin2010JCP, Xiu2004662, Acharjee20061934}) to address the problem posed by large $d$. In particular, we focus on partial differential equations (PDEs) with spatially dependent random coefficients. If considered in the whole spatial domain, the randomness is very high-dimensional. However, a few dominant random parameters could suffice locally~\cite{Chen2015}. As a result, we can obtain low-dimensional {\it local} representations of the random solution in each subdomain. To obtain the low dimensional representation in each subdomain, we use the Hilbert space KL expansion \cite{RK:Doostan2007}. Then, we reconstruct a global solution by stitching together the local solutions for each subdomain. 

The paper is organized as follows: Section \ref{RK:sec:spde} presents the problem of uncertainty quantification for PDEs with random coefficients. Section \ref{RK:sec:stochreduct} introduces our approach that combines the  domain decomposition and basis adaptation. Section \ref{RK:sec:numerical} contains numerical results for a two-dimensional steady-state diffusion equation for two different types of boundary conditions. Section \ref{conclusions} presents some conclusions and ideas for future work.

\section{PDEs with random coefficients}
\label{RK:sec:spde}

Let  $D$ be an open subset of $\mathbb{R}^n$ and $(\Omega, \Sigma, P)$ be a complete probability space 
with sample space $\Omega$, $\sigma$-algebra $\Sigma$, and  probability measure $P.$ We want to
find a random field, $u(x,\omega):D \times \Omega \rightarrow \mathbb{R},$ such that $P$-almost surely
%\begin{equation}\label{RK:eq:spdeop}
% \mathcal{L}(x,u;a(x,\omega)) = f(x,\omega)  \;\; \rm{in}~D\times \Omega,
%\end{equation}
%subject to the boundary condition
%\begin{equation}\label{RK:eq:spdebc}
% \mathcal{B}(x,u;a(x,\omega)) = h(x,\omega)  \;\; \rm{on}~\partial D\times \Omega,
%\end{equation}
\begin{align}\label{RK:eq:spdeop}
 \mathcal{L}(x,u(x,\omega);a(x,\omega)) &= f(x,\omega)  \;\; \rm{in}~D\times \Omega, \nonumber \\
 \mathcal{B}(x,u(x,\omega);a(x,\omega)) &= h(x,\omega)  \;\; \rm{on}~\partial D\times \Omega,
\end{align}
where $\mathcal{L}$ is a differential operator and $\mathcal{B}$ is a boundary operator. We model the uncertainty in the stochastic PDE (\ref{RK:eq:spdeop}) by treating the coefficient $a(x,\omega)$ as a random field and compute the effect of this uncertainty on the solution field $u(x,\omega).$ To solve the stochastic PDE numerically, we discretize the random fields $a(x,\omega)$ and $u(x,\omega)$ in both spatial and stochastic domains. In this paper, we will focus on the special case where $a(x,\omega)$ is a log-normal random field~\cite{Ghanem1999}. This means that $a(x,\omega) = \exp[g(x,\omega)]$, where $g(x,\omega)$ is a Gaussian random field whose mean and covariance function are known. We approximate $g(x,\omega)$ with a truncated KL expansion, while the coefficient $a(x,\omega)$ and solution $u(x,\omega)$ are approximated through truncated PC expansions. 

\subsection{Karhunen--Lo\`eve expansion of the random field $g(x,\omega)$}\label{RK:sec:kle}
The random field $g(x,\omega)$ can be approximated with a truncated KL expansion~\cite{RK:Loeve1977},
 \begin{equation}\label{RK:eq:gkl}
 g(x,\omega) \approx g(x, \boldsymbol{\xi}(\omega)) = g_0(x) + \sum_{i=1}^d \sqrt{\lambda_i} g_i(x) \xi_i(\omega),
\end{equation}
where $d$ is the number of random variables in the truncated expansion; $\boldsymbol{\xi} = (\xi_1, \cdots, \xi_d)^T$,
 $\xi_i$  are uncorrelated random variables with zero mean;
 $g_0(x)$ is the mean of the random field $g(x,\omega)$; and $(\lambda_i, g_i(x))$ are eigenvalues and eigenvectors obtained by solving the eigenvalue problem,
 \begin{equation}\label{RK:eq:eig}
 	\int_D C_g(x_1,x_2)  g_i(x_2)dx_2 = \lambda_i g_i(x_1),
\end{equation}
where $C_g(x_1,x_2)$ is the covariance function of the Gaussian random field $g(x,\omega)$. The eigenvalues are positive and non-increasing, and the eigenfunctions $g_i(x)$ are orthonormal, 
 \begin{equation}\label{RK:eq:ortho}
 	\int_D g_i(x) g_j(x)dx = \delta_{ij},
\end{equation}
where $\delta_{ij}$ is the Kronecker delta. In this work, we assume the random variables $\xi_i$ have Gaussian distribution. Therefore, $\xi_i$ are independent.

\subsection{Polynomial chaos expansion}\label{RK:sec:pce}
We approximate the input random field, $a(x,\omega)$ and the solution field $u(x,\omega)$ using truncated PC expansions~\cite{Cameron1947} in Gaussian random variables as follows:

\begin{equation}\label{RK:eq:pce_a}
 a(x,\omega) \approx  a(x, \boldsymbol{\xi}(\omega)) = a_0(x) + \sum_{i=1}^{N_{\xi}} a_i(x) \psi_i(\boldsymbol{\xi})
\end{equation}

and

\begin{equation}\label{RK:eq:pce_u}
  u(x,\omega) \approx u(x,\boldsymbol{\xi}(\omega)) \approx u_0(x) + \sum_{i=1}^{N_{\xi}} u_i(x) \psi_i(\boldsymbol{\xi}),
\end{equation}
where  $N_{\xi} = \frac{(d+p)!}{d!~p!}$ is the number of terms in PC expansion for dimension $d$ and order $p$, $u_0(x)$ is the mean of the solution field, $u_i(x)$ are PC coefficients, and $\{\psi_i(\boldsymbol{\xi})\}$ are multivariate Hermite polynomials. These polynomials are orthogonal with respect to the inner product defined by the expectation in the stochastic space, 

\begin{equation}\label{RK:eq:innprod}
\langle \psi_i(\boldsymbol{\xi}), \psi_j(\boldsymbol{\xi}) \rangle \equiv \int_{\Omega} \psi_i(\boldsymbol{\xi}(\omega)) \psi_j(\boldsymbol{\xi}(\omega)) dP(\omega) = \delta_{ij},
\end{equation}
where $dP(\omega)$ is the diagonal $d$-dimensional Gaussian distribution. Once the random coefficient $a(x,\omega)$ and the solution field $u(x,\omega)$ are approximated, the stochastic PDE \eqref{RK:eq:spdeop} is transformed into 
\begin{align}\label{RK:eq:spdexi}
 \mathcal{L}(x,\boldsymbol{\xi}, u(x,\boldsymbol{\xi});a(x,\boldsymbol{\xi})) &= f(x,\boldsymbol{\xi})  \;\; \rm{in}~D\times \Omega, \nonumber \\
 \mathcal{B}(x,\boldsymbol{\xi}, u(x,\boldsymbol{\xi});a(x,\boldsymbol{\xi})) &= h(x,\boldsymbol{\xi})  \;\; \rm{on}~\partial D\times \Omega.
\end{align}
Intrusive methods, such as stochastic Galerkin~\cite{RK:Ghanem1991}, or non-intrusive methods, such, as sparse-grid collocation~\cite{RK:Xiu2002}, can be used to solve the parameterized stochastic PDE \eqref{RK:eq:spdexi}. However, both methods suffer from the curse of dimensionality, i.e., their computational costs  exponentially increase with  increasing $d$ and/or  $p$. To address the curse of dimensionality, we  propose a novel approach based on the basis adaptation method, introduced in~\cite{RK:Tipireddy2014} to compute locally accurate solutions. Our approach combines  the basis adaptation method with spatial domain decomposition, which allows us to compute an accurate solution in the entire domain. In the following section, we demonstrate how basis adaptation can be applied to find a low-dimensional solution in non-overlapping spatial subdomains.

\section{Basis adaptation in a spatial subdomain} 
\label{RK:sec:stochreduct}
To address the curse of dimensionality in solving~\eqref{RK:eq:spdeop}, we first decompose the spatial domain $D$ into a set of non-overlapping subdomains. Then, we use the basis adaptation to find a low-dimensional solution space in each subdomain. 
\subsection{Domain decomposition}
\label{RK:sec:dd}
We divide the spatial domain $D \subset \mathbb{R}^n$ into a set of non-overlapping subdomains~\cite{Chen2015}, $D_{s} \subset D, s = 1, \cdots, S,$ such that 
\begin{equation}\label{RK:eq:domaind}
	{D} = \bigcup_{s=1}^S {D_s}, \quad D_s \cap D_{s^{\prime}}=\emptyset, s \neq s^{\prime}. 
\end{equation}	
In each subdomain $D_s$, we find a low-dimensional stochastic basis (PC) for a new set of random variables $\tilde{\boldsymbol{\eta}}^s = \{\eta_1^s, \cdots, \eta_{r}^s \}, r \ll d,$ adapted to that subdomain, such that the solution can be well represented in the low-dimensional basis. When the stochastic basis is adapted to the subdomain $D_s$, the low-dimensional solution computed by solving~\eqref{RK:eq:spdeop} is accurate only in $D_s$. Hence, we repeat the process in each subdomain to get an accurate low-dimensional solution in the entire spatial domain.  

\subsection{Basis adaptation through Hilbert-Karhunen-Lo\`eve expansion}
\label{RK:sec:hkle}
Here, we extend the basis adaptation method of~\cite{RK:Tipireddy2014} to construct a low-dimensional solution representation in each subdomain $D_s$ using the Hilbert-Karhunen-Lo\`eve expansion.
The KL expansion~\cite{RK:Loeve1977} described in Section~\ref{RK:sec:kle} provides optimal representation of a random field in $L_2$ space. However, to take advantage of this property for PDEs with random coefficients, the solution should satisfy certain regularity and smoothness conditions~\cite{RK:Doostan2007}. Thus, it is appropriate to seek a representation of the solution field $u(x,\omega)$ in a subset of $L_2(\Omega)$. Such a restricted expansion has been explored in various applications~\cite{RK:Levy1999, RK:Kirby1992, RK:Silverman1996, Berkooz1993539, Christensen1999}. 

In this work, we first compute the Gaussian part (up to linear terms in PC expansion) of the solution in the whole domain $D,$ i.e., 
\begin{equation}\label{RK:eq:pce_ug2}
	 u_g(x,\boldsymbol{\xi}(\omega)) = u_0(x) + \sum_{i=1}^{d} u_i(x) \xi_i.
\end{equation}
Note that in realistic situations, this computation could be feasible even if the full solution computation is not. 

Let $u^s_g(x,\boldsymbol{\xi}(\omega)) \subset u_g(x,\boldsymbol{\xi}(\omega))$ be the Gaussian part of the solution in subdomain $D_s$,  i.e.,
\begin{equation}\label{RK:eq:pce_ug_s}
	u^s_g(x,\boldsymbol{\xi}(\omega)) = u_g(x,\boldsymbol{\xi}(\omega)) \mathbb{I}_{D_s}(x), 
\end{equation}
where $\mathbb{I}$ is the indicator function such that for any set $\alpha$ ($\mathbb{I}_{\alpha} = 1$, if $x \in \alpha,$ and $\mathbb{I}_{\alpha} = 0$, if $x \notin \alpha$). For each subdomain $D_s,$ we construct the covariance function of the Gaussian part of the solution $u^s_g(x,\boldsymbol{\xi})$ as follows: 
\begin{equation}\label{RK:eq:cov_ug}
	 C^s_{u_g}(x_1,x_2) = \sum_{i=1}^{d} u_i(x_1) u_i(x_2), \quad x_1, x_2 \in\bar{D_s}.
\end{equation}
The Hilbert space KL expansion of $u^s_g(x,\boldsymbol{\xi})$ (see e.g., \cite{RK:Doostan2007}) is
\begin{equation}\label{RK:eq:ugkl}
	u^s_g(x,\boldsymbol{\xi}(\omega)) = u^s_0(x) + \sum_{i=1}^{d} \sqrt{\mu^s_i} \phi^s_i(x) \eta^s_i(\omega), \quad x \in \bar{D_s}
\end{equation}
where $(\mu^s_i, \phi^s_i(x))$ are eigenpairs corresponding to the above expansion. They can be obtained by solving the eigenvalue problem: 
\begin{equation}\label{RK:eq:cov_cg}
	 \int_{D_s} C^s_{u_g}(x_1,x_2) \phi^s_i(x_1)dx_1 = \mu^s_i \phi^s_i(x_2), \quad i = 1, 2, \cdots. 
\end{equation}
After obtaining the solution of the eigenvalue problem, we can use \eqref{RK:eq:ugkl} to write the random variable $\eta^s_i$ as
\begin{align}\label{RK:eq:pce_eta}
	 \eta^s_i &= \frac{1}{\sqrt{\mu^s_i}}\int_{D_s} (u^s_g(x,\boldsymbol{\xi}) - \tilde{u}^s_0(x)) \phi^s_i(x)dx, \nonumber \\
	 		&= \frac{1}{\sqrt{\mu^s_i}}\int_{D_s} (u_0(x) + \sum_{j=1}^{d} u_j(x) \xi_j - \tilde{u}^s_0(x)) \phi^s_i(x)dx. \quad i = 1, 2, \cdots, d
\end{align}
Note that $u_0(x) = u^s_0(x)$ by construction. Hence, 
 \begin{align}\label{RK:eq:eta_xi}
	 \eta^s_i &= \frac{1}{\sqrt{\mu^s_i}}\int_{D_s} \left ( \sum_{j=1}^{d} u_j(x) \xi_j \right) \phi^s_i(x)dx, \quad x \in D_s, i = 1, 2, \cdots, d \nonumber \\
	 		&= \sum_{j=1}^{d} \left ( \frac{1}{\sqrt{\mu^s_i}}\int_{D_s}  u_j(x) \phi^s_i(x)dx  \right)\xi_j, \quad x \in D_s, i = 1, 2, \cdots, d \nonumber \\
	 		&= a^s_{ij} \xi_j,
\end{align}
where $a^s_{ij} =  \frac{1}{\sqrt{\mu^s_i}}\int_{D_s}  u_j(x) \phi^s_i(x)dx, i,j = 1, \cdots, d.$ This provides a linear relation between $\eta^s_i$ and $\{\xi_j\}$. Because, $\xi_j$ are standard Gaussian random variables, the new variables $\{\eta^s_i\}$ also are Gaussian random variables and can be normalized to obtain standard Gaussian random variables. The resulting Gaussian random variables $\{\eta^s_i\}$ are uncorrelated and, thus, also independent. We can reformulate the stochastic PDE~\eqref{RK:eq:spdeop} in terms of $\{\eta^s_i\}$ and solve using intrusive or non-intrusive methods.

\subsection{Dimension reduction in the subdomain $D_s \subset D$}
\label{RK:sec:ba}
Let $A_s = [a^s_{ij}]$ be an isometry in $\mathbb{R}^d$ and define $\boldsymbol{\eta}^s$ as
\begin{equation}\label{RK:eq:eta}
	\boldsymbol{\eta}^s =A_s \boldsymbol{\xi}, \quad A_s {A_s}^T = \boldsymbol{I}, 
\end{equation}
where $\boldsymbol{\eta}^s = \{\eta^s_1, \cdots, \eta^s_d\}^T$ is a vector of standard normal random variables. With the mapping in \eqref{RK:eq:eta}, $\boldsymbol{\eta}^s$ is defined in the Gaussian Hilbert space spanned by random variables $\boldsymbol{\xi}$ in $D_s.$ This means
\begin{equation}\label{RK:eq:u_tilde}
	u(x,\boldsymbol{\xi}(\omega)) = {u}^{A_s}(x,\boldsymbol{\eta}^s(\omega))
	= {u}^{A_s}_0(x) + \sum_{i=1}^{N_{\eta^s}} {u}^{A_s}_i(x) \psi_i(\boldsymbol{\eta}^s), \quad x \in D_s.
\end{equation}

For this construction to be worth the effort, we should be able to obtain an accurate representation of the solution without having to use functions of all $d$ new variables $\eta^s_i, \; i=1,\ldots,d.$ This depends on how fast the magnitude of the eigenvalues $\mu^s_i$ ($i=1,\ldots,d$) decays. For the steady-state diffusion equation considered here, the solution is expected to be smoother than the input random field $a(x,\omega).$ In addition, we employ the Hilbert space KL expansion in a spatial subdomain $D_s \subset D$. Thus, eigenvalues $\mu^s_i$ in subdomain $D_s$ are expected to decay faster than the eigenvalues $\lambda_i$ of the input random field (see Fig. \ref{RK:fig:eigen}). Because of this, it is sufficient to perform the Hilbert space KL expansion in \eqref{RK:eq:ugkl} only in $r$($\ll d$) dimensions. 

With the new set of random variables $\tilde{\boldsymbol{\eta}}^s = \{\eta_1, \cdots, \eta_r \}$, \eqref{RK:eq:spdexi} can be reformulated as 
\begin{align}\label{RK:eq:spdeeta}
 \mathcal{L}(x,\boldsymbol{\xi}, \tilde{\boldsymbol{\eta}}^s, u(x,\tilde{\boldsymbol{\eta}}^s);a(x,\boldsymbol{\xi})) &= f(x,\boldsymbol{\xi})  \;\; \rm{in}~D\times \Omega, \nonumber \\
 \mathcal{B}(x,\boldsymbol{\xi}, \tilde{\boldsymbol{\eta}}^s, u(x,\tilde{\boldsymbol{\eta}}^s);a(x,\boldsymbol{\xi})) &= h(x,\boldsymbol{\xi})  \;\; \rm{on}~\partial D\times \Omega.
\end{align}
The stochastic PDE \eqref{RK:eq:spdeeta} can be solved using non-intrusive or intrusive methods (described in the two sections that follow).

\subsection{Sub-domain solution using intrusive methods}
Expanding the solution of \eqref{RK:eq:spdeeta} with Hermite polynomials in $\tilde{\boldsymbol{\eta}}^s$ and the random coefficient $a(x,\boldsymbol{\xi})$ with Hermite polynomials in $\boldsymbol{\xi}$, and performing a Galerkin projection leads to the deterministic PDEs,
\begin{equation}\label{RK:eq:sg_eta}
	 \sum_{j=0}^{N_{\eta^s}} \sum_{i=0}^{N_{\xi}}  \mathcal{L}_{i}[\tilde{u}^{A_s}_{j}] E[\psi_{i}(\boldsymbol{\xi}) \psi_{j}(\tilde{\boldsymbol{\eta}}^s)  \psi_{k}(\tilde{\boldsymbol{\eta}}^s)] = 0, \quad  k = 1, \cdots, N_{\eta^s},
\end{equation}
for the expansion coefficients $\tilde{u}^{A_s}_{j}, \; j=1,\ldots,N_{\eta^s}.$ Also, in \eqref{RK:eq:sg_eta}, $\mathcal{L}_{i} = (\mathcal{L}, \psi_i(\boldsymbol{\xi}) )$ is the inner product of the operator $\mathcal{L}$ with the PC coefficient $\psi_i(\boldsymbol{\xi})$.
If the number of retained random variables $\eta^s_i$ is small, then the number of PC terms  $N_{\eta^s}$ is much smaller than $N_{\xi}.$ 

Although we find a low-dimensional stochastic basis adapted to a subdomain $D_s$, we solve \eqref{RK:eq:sg_eta} in the entire spatial domain $D$. However, because the random variables $\eta^s_i$ are adapted to the subdomain $D_{s},$ the solution is accurate only in that subdomain. Hence, for each subdomain $D_s, \; s \in S$, we solve \eqref{RK:eq:sg_eta} with stochastic basis adapted to that subdomain and retain the solution corresponding to that subdomain while  discarding the solution for the remaining subdomains. In this way, we obtain a solution for each subdomain using the adaptive basis in that subdomain. 

Let $|S|$ be the total number of non-overlapping subdomains. If the subdomains are of equal size (and resolution), the computational cost for obtaining the solution in the entire domain is $|S|$ times the cost of solving \eqref{RK:eq:sg_eta} plus the cost of solving \eqref{RK:eq:spdexi} up to the first order PC.  

\subsection{Sub-domain solution using non-intrusive methods}
In a non-intrusive method, a PDE is solved at predefined quadrature points (in random space). We use sparse-grid collocation points based on Smolyak approximation~\cite{RK:Smolyak1963, RK:Nobile2008}. Unlike the tensor-product of one dimensional quadrature points, the sparse-grid method judiciously chooses products with only a small number of quadrature points. These product rules depend on an integer value called sparse-grid level~(see~\cite{RK:Smolyak1963, RK:Nobile2008}). As the order of PC expansion increases, we need higher sparse-grid levels to maintain solution accuracy. In the sparse-grid method, the number of collocation points increases with the stochastic dimension (number of random variables) and  the sparse-grid level. Let $\tilde{\boldsymbol{\eta}}_q^s$ be the collocation points associated with the random variables $\tilde{\boldsymbol{\eta}}^s$. Then, the corresponding points ${\boldsymbol{\xi}}_q^s$ are obtained as 
\begin{equation}\label{RK:eq:xi_Aeta}
	{\boldsymbol{\xi}}^s_q =   [A_s^{-1}]_{r} \tilde{\boldsymbol{\eta}}^s_q,
\end{equation}
where the $d \times r$ matrix $[A_s^{-1}]_r$ consists of the first $r$ columns of $A_s^{-1}$. Finally, the deterministic PDE 
\begin{equation}\label{RK:eq:Lu_eta}
	\mathcal{L} [{u}({{\boldsymbol{\xi}}_q^s})] =   0	
\end{equation}
is solved for each collocation point ${{\boldsymbol{\xi}}_q^s}$ ($q = 1,\cdots, Q$), and the PC coefficient $\tilde{u}^{A_s}_i$ is computed by projection,
\begin{equation}\label{RK:eq:ui_eta}
	\tilde{u}^{A_s}_i(x) = \sum_{q=1}^Q {u}({\boldsymbol{\xi}}^s_q) \psi_i(\tilde{\boldsymbol{\eta}}^s_q) w_q,
\end{equation}
where $w_q$ are weights for the quadrature points. 

\subsection{The error of the low-dimensional solution, $\tilde{u}^{A_s}(x,\tilde{\boldsymbol{\eta}}^s)$} 
\label{RK:sec:baerror}
In the basis adaptation method, if the dimension of the new basis in $\boldsymbol{\eta}^s$ and the dimension of original basis in $\boldsymbol{\xi}$ are equal, there is no dimension reduction, and the accuracy of the solution $\tilde{u}^{A_s}(x,\boldsymbol{\eta}^s)$ in terms of $\boldsymbol{\eta}^s$ is the same as that of the solution $u(x,\boldsymbol{\xi}).$ This means: 
\begin{align}\label{RK:eq:error}
	\epsilon(x, \omega) = u(x,\boldsymbol{\xi}(\omega)) - \tilde{u}^{A_s}(x,\boldsymbol{\eta}^s(\omega)) = 0, \quad \text{if } r=d.
\end{align} 
However, if the eigenvalues $\mu^s_i$ decay fast (see Fig.~\ref{RK:fig:eigen}), we can truncate the Hilbert space KL expansion in \eqref{RK:eq:ugkl} to $r+1$ terms and represent the solution in a low-dimensional ($r \ll d$) space, 
\begin{align}\label{RK:eq:error2}
	u(x,\boldsymbol{\xi}(\omega)) &= \tilde{u}^{A_s}_0(x) + \sum_{\mathcal{I}_1} \tilde{u}_{\mathcal{I}_1} \psi_{\mathcal{I}_1}(\tilde{\boldsymbol{\eta}}^s) + \sum_{\mathcal{I}_2} \tilde{u}_{\mathcal{I}_2} \psi_{\mathcal{I}_2}(\hat{\boldsymbol{\eta}}^s) \nonumber \\
	u(x,\boldsymbol{\xi}(\omega)) &- \left ( \tilde{u}^{A_s}_0(x) + \sum_{\mathcal{I}_1} \tilde{u}_{\mathcal{I}_1} \psi_{\mathcal{I}_1}(\tilde{\boldsymbol{\eta}}^s)\right) = \sum_{\mathcal{I}_2} \tilde{u}_{\mathcal{I}_2} \psi_{\mathcal{I}_2}(\hat{\boldsymbol{\eta}}^s) \nonumber \\	
	u(x,\boldsymbol{\xi}(\omega)) &- \tilde{u}^{A_s}(x,\tilde{\boldsymbol{\eta}}^s(\omega)) = \sum_{\mathcal{I}_2} \tilde{u}_{\mathcal{I}_2} \psi_{\mathcal{I}_2}(\hat{\boldsymbol{\eta}}^s) \nonumber \\
	\epsilon(x, \omega) &= \sum_{\mathcal{I}_2} \tilde{u}_{\mathcal{I}_2} \psi_{\mathcal{I}_2}(\hat{\boldsymbol{\eta}}^s), \quad \text{if } r \ll d
\end{align} 
where, $\tilde{\boldsymbol{\eta}}^s = \{\eta^s_1, \cdots,\eta^s_r \}$, $\hat{\boldsymbol{\eta}}^s = \{\eta^s_{r+1}, \cdots,\eta^s_d \}$, and multi-indices, $\mathcal{I}_1$ and $\mathcal{I}_2,$ correspond to PC expansion terms in $\tilde{\boldsymbol{\eta}}^s$ and $\hat{\boldsymbol{\eta}}^s$, respectively.

\section{Numerical examples}
\label{RK:sec:numerical}
In this section, we apply the proposed approach to a two-dimensional steady-state diffusion equation with random coefficient, defined on the spatial domain $D = [0, 240] \times [0, 60]$ such that $x=(x_1,x_2) \in D$. We solve the equation for two types of boundary conditions. In the first case, we use Dirichlet boundary conditions at the boundaries perpendicular to the  $x_1$ direction and Neumann boundary conditions at the other two boundaries. In the second case, we use Dirichlet boundary conditions on all four sides. 

\begin{figure}[!h]
\begin{center}
{\includegraphics[scale=.6]{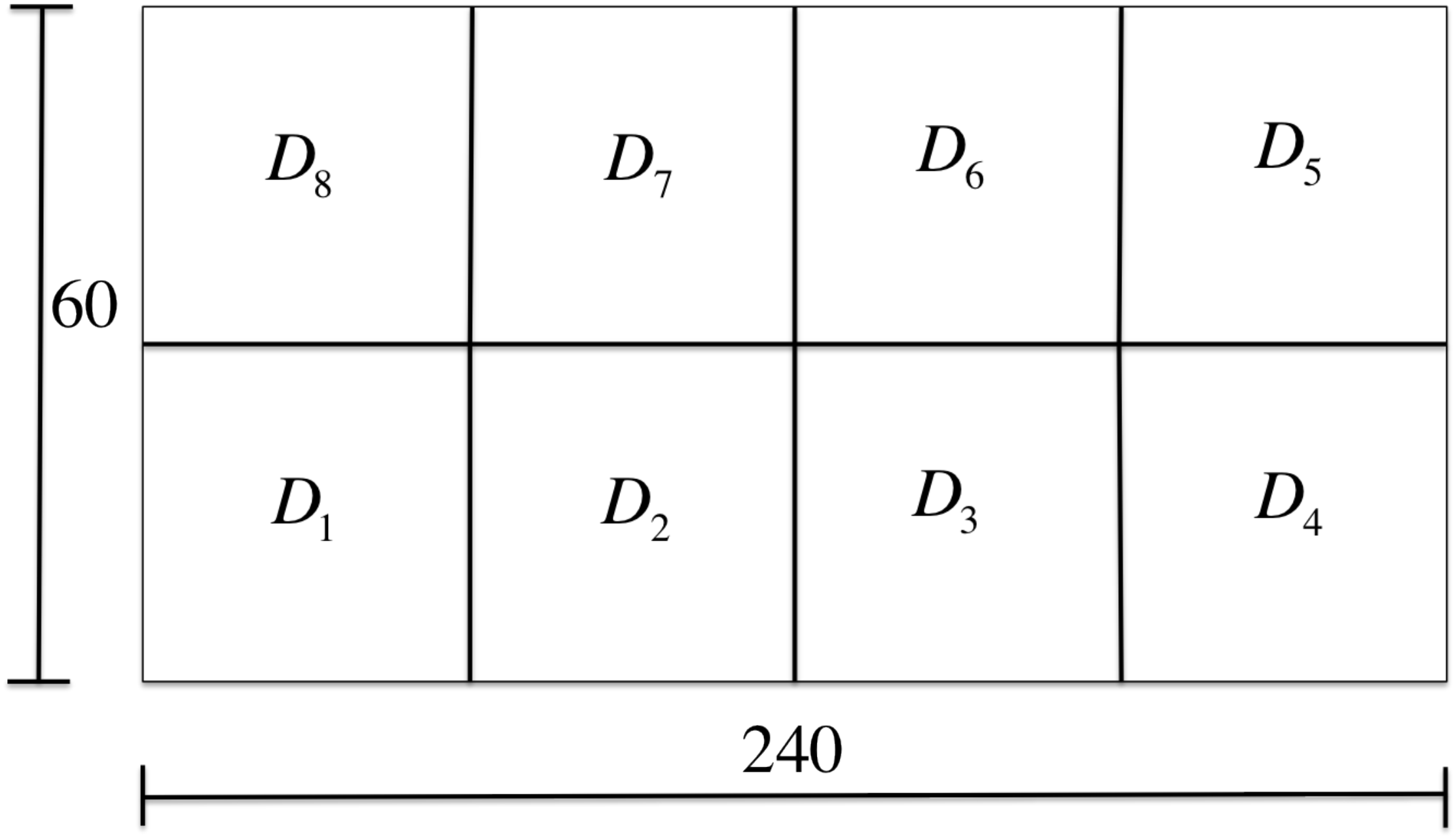} }
 \caption{Spatial domain decomposed into eight non-overlapping subdomains.} \label{RK:fig:DD}
\end{center}
\end{figure}

\begin{figure}[!h]
\begin{center}
{\includegraphics[scale=.24]{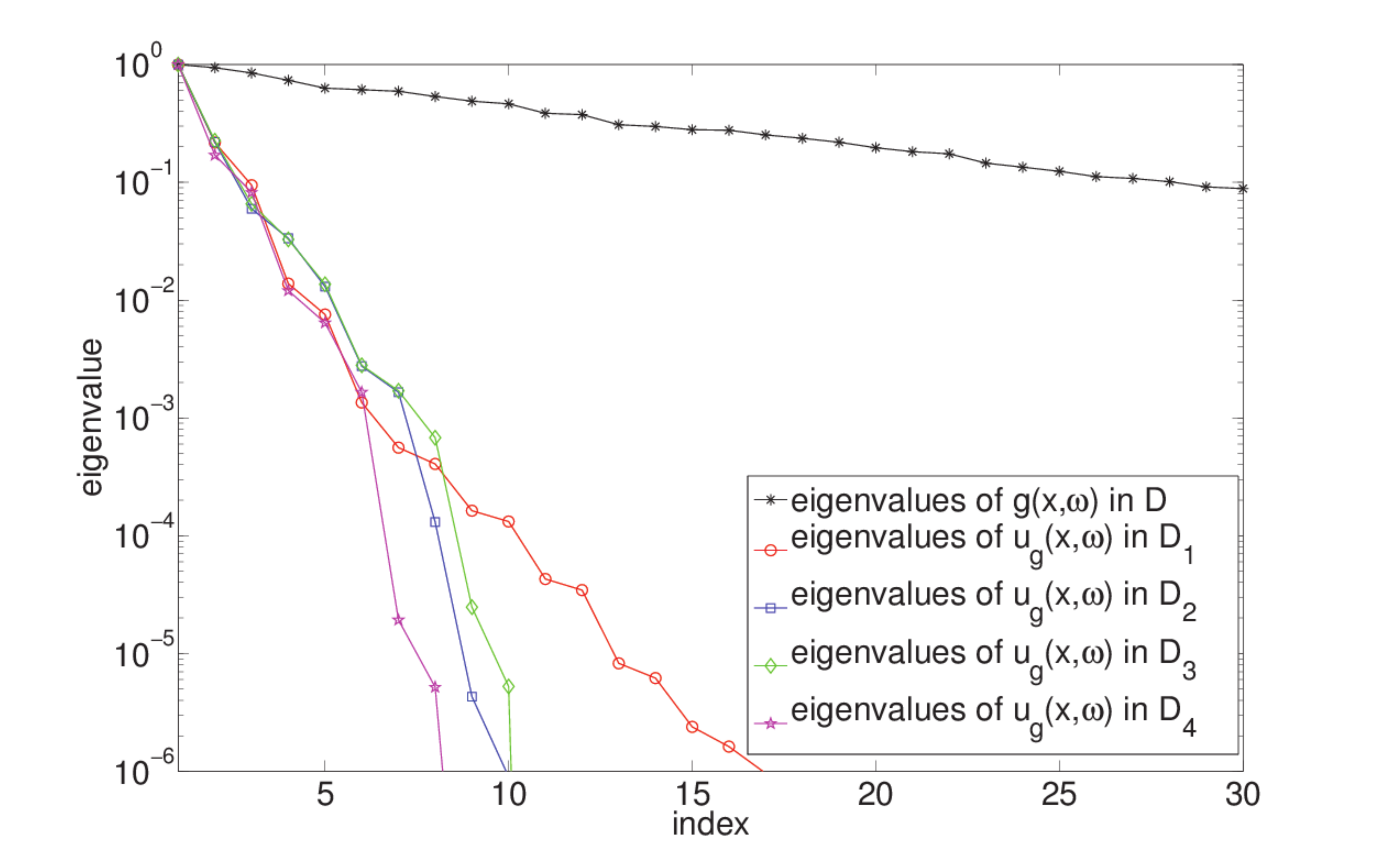} }
 \caption{Decay of the eigenvalues of covariance function of $g(x,\omega)$ in domain $D$ and the covariance function of Gaussian solution from subdomains $D_1, D_2, D_3,$ and $D_4.$} \label{RK:fig:eigen}
\end{center}
\end{figure}

Let the random coefficient $a(x,\omega):D \times \Omega \rightarrow \mathbb{R}$ be bounded and strictly positive, 
\begin{equation}\label{RK:eq:rf_bound}
	0 < a_l \leq a(x,\omega) \leq a_u < \infty \quad \rm{a.e.} \quad \rm{in} \quad D\times \Omega.
\end{equation}
In the first case, we solve the boundary-value problem:
\begin{align}\label{RK:eq:spde}
	-\nabla . (a(x,\omega) \nabla u(x,\omega))&=f(x,\omega) \;\; \rm{in}~D\times \Omega,  \nonumber \\  
	u(x,\omega)&=100 \;\; \text{on}~x_1 = 0, \nonumber \\
	u(x,\omega)&=10 \;\; \text{on}~x_1 = 240, \nonumber \\
	\vec{n} \cdot \nabla u(x,\omega) & = 0  \;\; \text{on}~x_2 = 0, \nonumber \\
	\vec{n} \cdot \nabla u(x,\omega) & = 0  \;\; \text{on}~x_2 = 60,
\end{align} 
where  $u(x,\omega):D \times \Omega \rightarrow \mathbb{R}$.
In the second case, we solve the boundary-value problem:
\begin{align}\label{RK:eq:spde_dbc}
	-\nabla . (a(x,\omega) \nabla u(x,\omega))&=f(x,\omega) \;\; \rm{in}~D\times \Omega,  \nonumber \\  
	u(x,\omega)&=100 \;\; \text{on}~x_1 = 0, \nonumber \\
	u(x,\omega)&=10 \;\; \text{on}~x_1 = 240, \nonumber \\
	u(x,\omega)&=0 \;\; \text{on}~x_2 = 0, \nonumber \\
	u(x,\omega)&=0 \;\; \text{on}~x_2 = 60. 
\end{align} 
We assume that random coefficient $a(x,\omega) = \exp[g(x,\omega)]$ has log-normal distribution with mean $a_0(x) = 5.0$ and standard deviation $\sigma_a = 2.5.$ The Gaussian random field $g(x,\omega)$ has correlation function 
\begin{equation}\label{RK:eq:cov_a}
	C_g(x,y) = \sigma_g^2 \exp \left(-\frac{(x_1-y_1)^2}{l_1^2} -\frac{(x_2-y_2)^2}{l_2^2} \right) \quad \rm{in} \quad D\times \Omega,
\end{equation}
where the standard deviation $\sigma_g = \sqrt{\ln \left(1+\frac{\sigma_a}{a_0(x)^2}\right)},$ the mean $g_0(x) = \ln \left (\frac{a_0(x)}{(\sqrt{1+\frac{\sigma_a}{a_0(x)^2}})} \right),$ and the correlation lengths $l_1 = 24$ and $l_2 = 20.$
We decompose the spatial domain into eight subdomains and independently find a low-dimensional solution in each subdomain  using the KL expansion and  basis adaptation methods described in Sections~\ref{RK:sec:hkle} and \ref{RK:sec:ba}. To evaluate the accuracy of the solution obtained with the reduced model, we solve the full stochastic system in the domain $D$ with stochastic dimension 10 and sparse-grid level 5. This requires $Q=8761$ collocation points.  In the domain decomposition with basis adaptation approach, we initially solve the full system in the entire domain $D$ with a coarse grid (sparse-grid level 3) and dimension 10, which required $Q=221$ collocation points. We use a coarser grid to obtain the Gaussian component of the solution and identify the low-dimensional space directions in each spatial subdomain. The low-dimensional representation computation requires 165 collocation points (dimension = 3 and sparse-grid level 5) per subdomain. As result, the total number the corresponding deterministic \eqref{RK:eq:Lu_eta} that must be solved in the domain decomposition with basis adaptation method is $8 \times 165 + 221 = 1541$, which is much smaller than  $Q=8761$ needed to solve the problem without basis adaptation. Figure~\ref{RK:fig:DD} depicts the decomposition of the spatial domain into eight subdomains. Figure~\ref{RK:fig:eigen} compares the decay of eigenvalues of the covariance function of $g(x,\omega)$ in domain $D$ and those of the covariance function of the Gaussian solution $u_g(x,\omega)$ in subdomains $D_1, D_2, D_3,$ and $D_4$. This figure demonstrates that eigenvalues decay fast in all subdomains, which justifies using a low-dimensional solution in each subdomain after basis adaptation. 

\begin{figure}[t!]
    \centering
    \begin{subfigure}[t]{0.3\textwidth}
        \centering
        \includegraphics[height=1.2in]{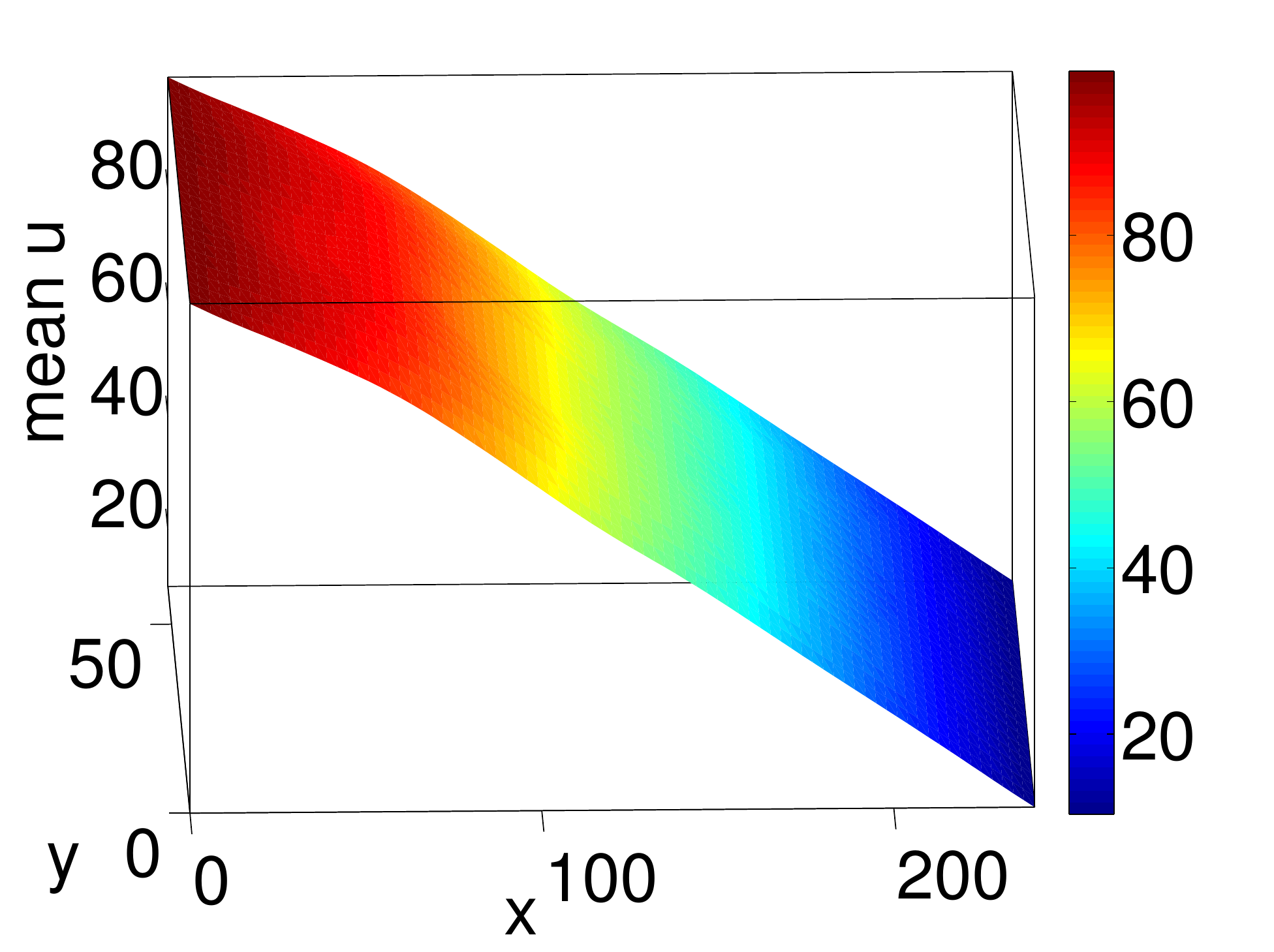}
        \caption{Mean, $\xi, d=10$} \label{RK:fig:u_mean_D1_full}
    \end{subfigure}        
    \begin{subfigure}[t]{0.3\textwidth}
        \centering
        \includegraphics[height=1.2in]{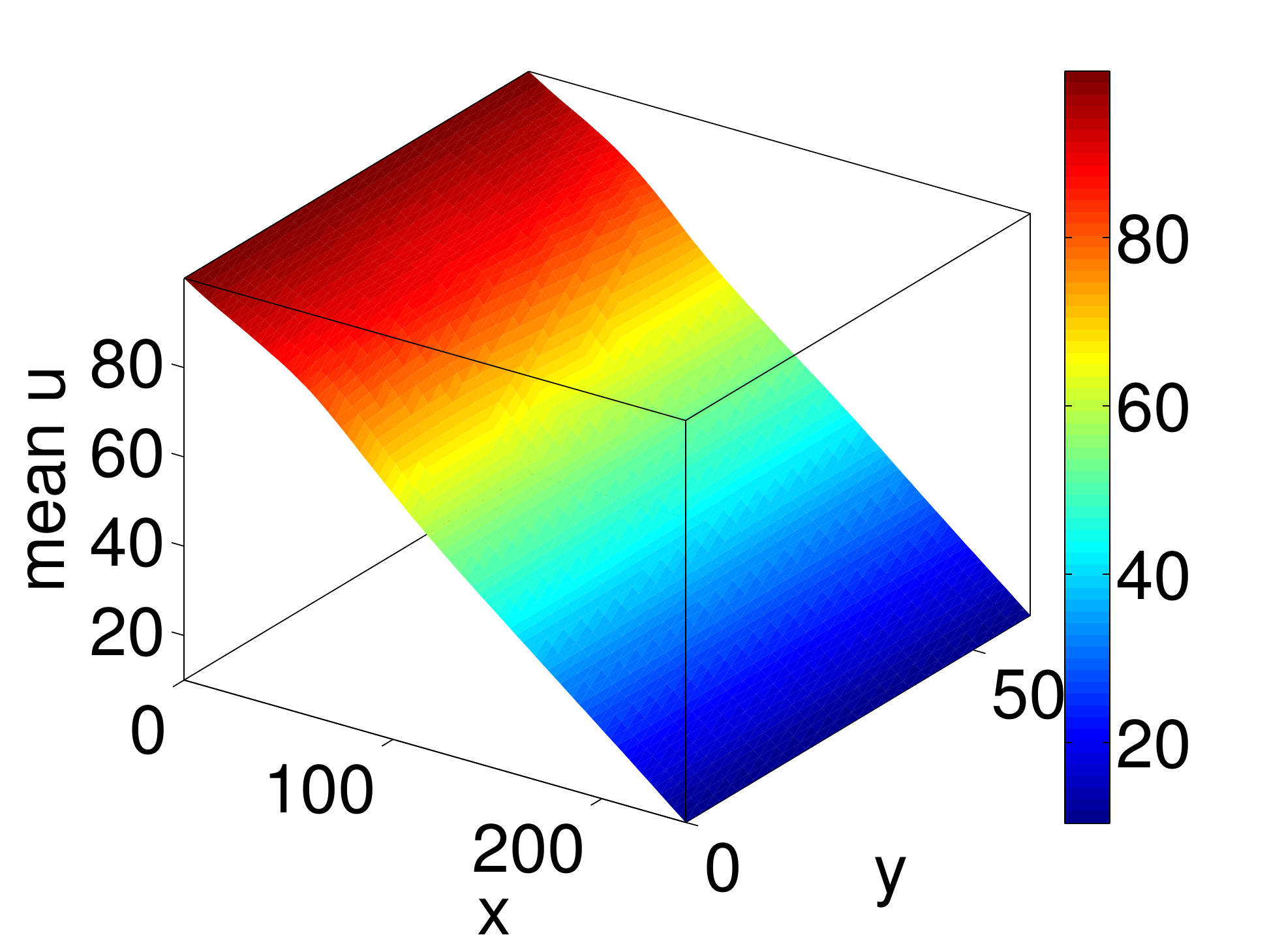}
        \caption{Mean, $\eta, d=3$} \label{RK:fig:u_mean_D1_reduced}
    \end{subfigure}    
    \begin{subfigure}[t]{0.3\textwidth}
        \centering
        \includegraphics[height=1.2in]{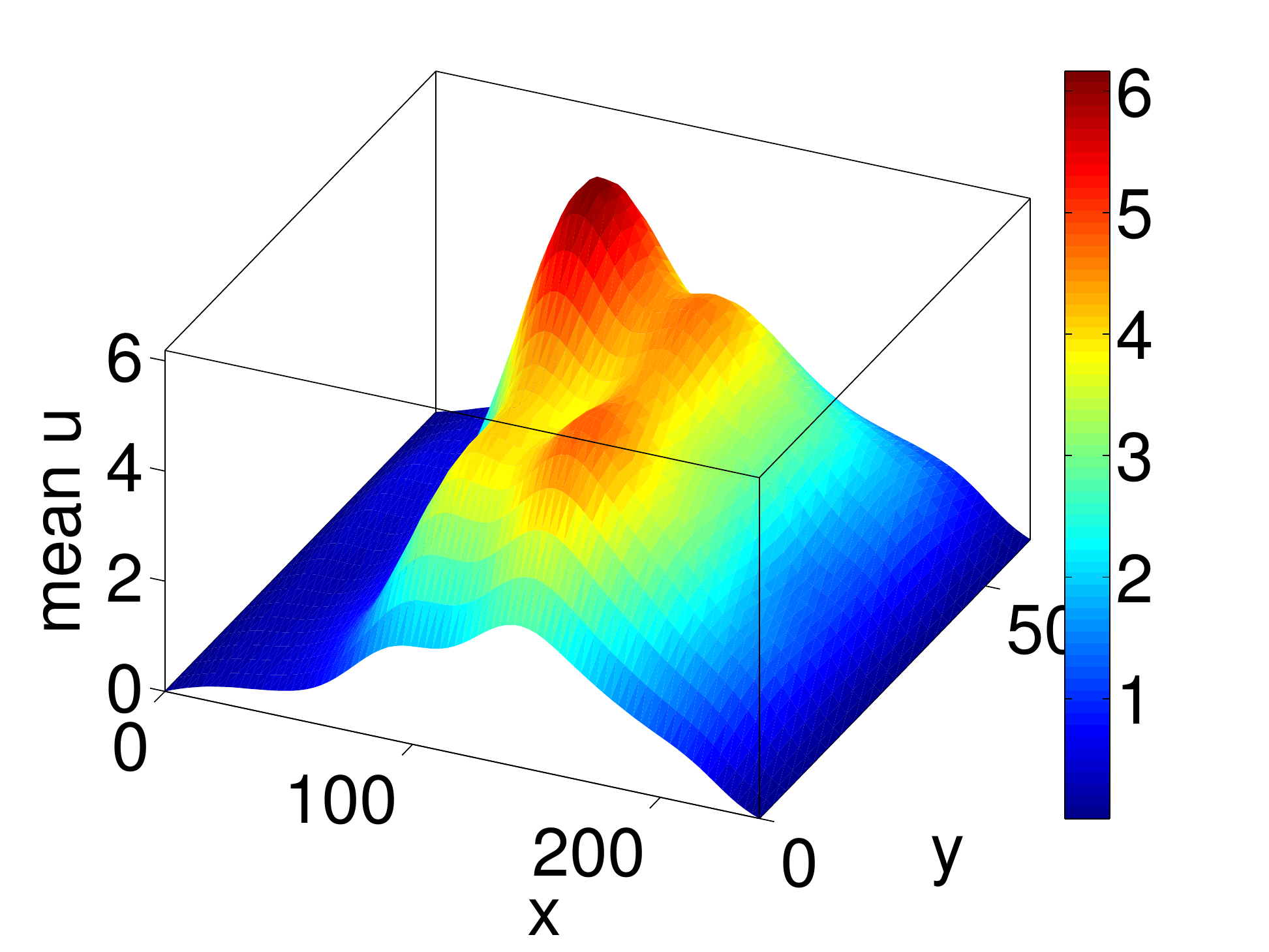}
        \caption{Error, $\eta, d=3$} \label{RK:fig:u_mean_D1_error}
    \end{subfigure}  
   
    \caption{Mean of the solution with the stochastic basis adapted for subdomain $D_1$ with random variables, $\eta$ and dimension = 3, order = 3.} \label{RK:fig:u_mean_D1}
\end{figure}

\begin{figure}[t!]
    \centering
    \begin{subfigure}[t]{0.3\textwidth}
        \centering
        \includegraphics[height=1.2in]{figures/u_mean_xid10_p3}
        \caption{Mean, $\xi, d=10$}\label{RK:fig:u_mean_D2_full}
    \end{subfigure}        
    \begin{subfigure}[t]{0.3\textwidth}
        \centering
        \includegraphics[height=1.2in]{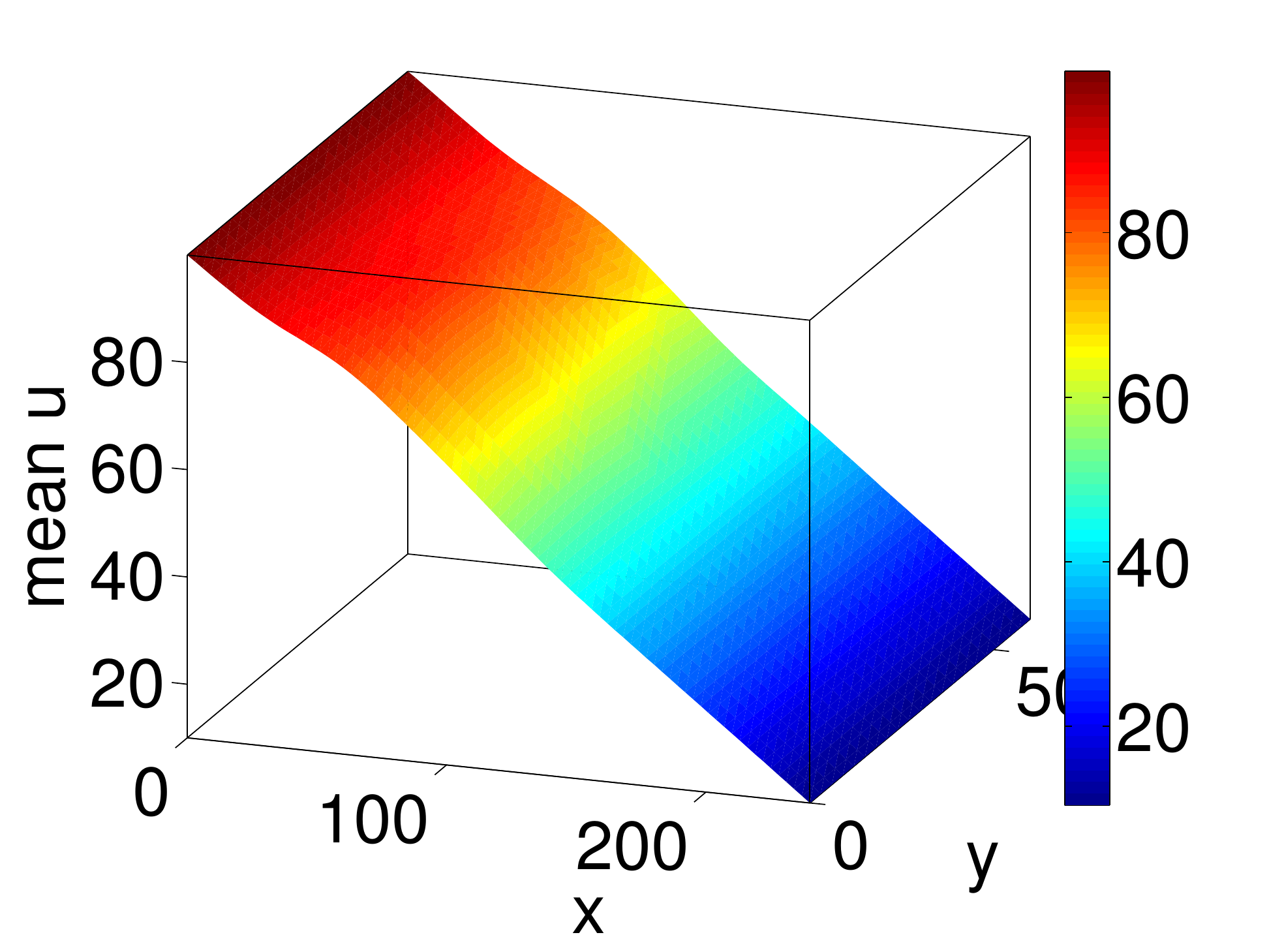}
        \caption{Mean, $\eta, d=3$}\label{RK:fig:u_mean_D2_reduced}
    \end{subfigure}    
    \begin{subfigure}[t]{0.3\textwidth}
        \centering
        \includegraphics[height=1.2in]{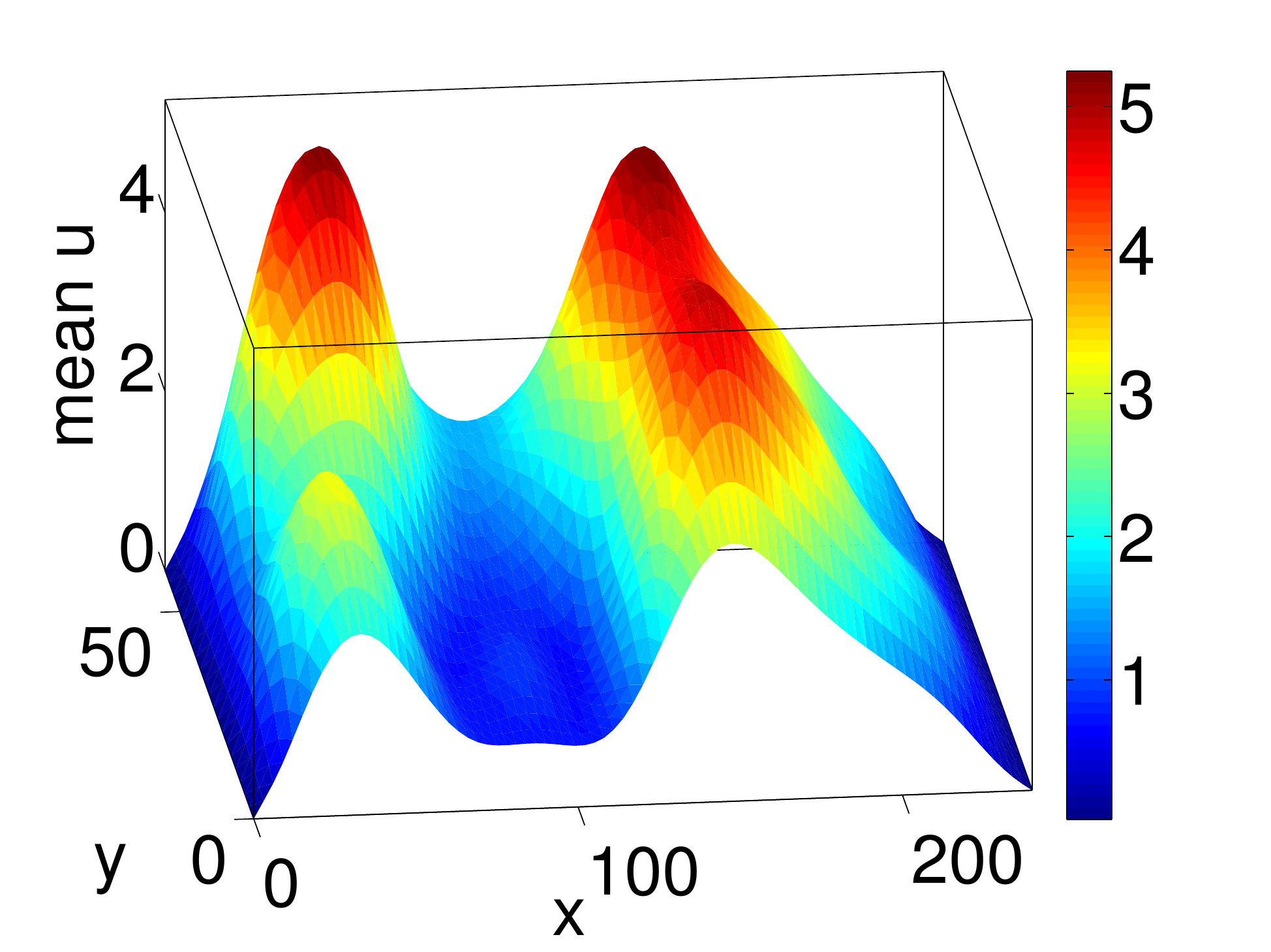}
        \caption{Error, $\eta, d=3$}\label{RK:fig:u_mean_D2_error}
    \end{subfigure}  

    \caption{Mean of the solution with the stochastic basis adapted for subdomain $D_2$ with random variables, $\eta$ and dimension = 3, order = 3.} \label{RK:fig:u_mean_D2}
\end{figure}

\begin{figure}[t!]
    \centering
    \begin{subfigure}[t]{0.3\textwidth}
        \centering
        \includegraphics[height=1.2in]{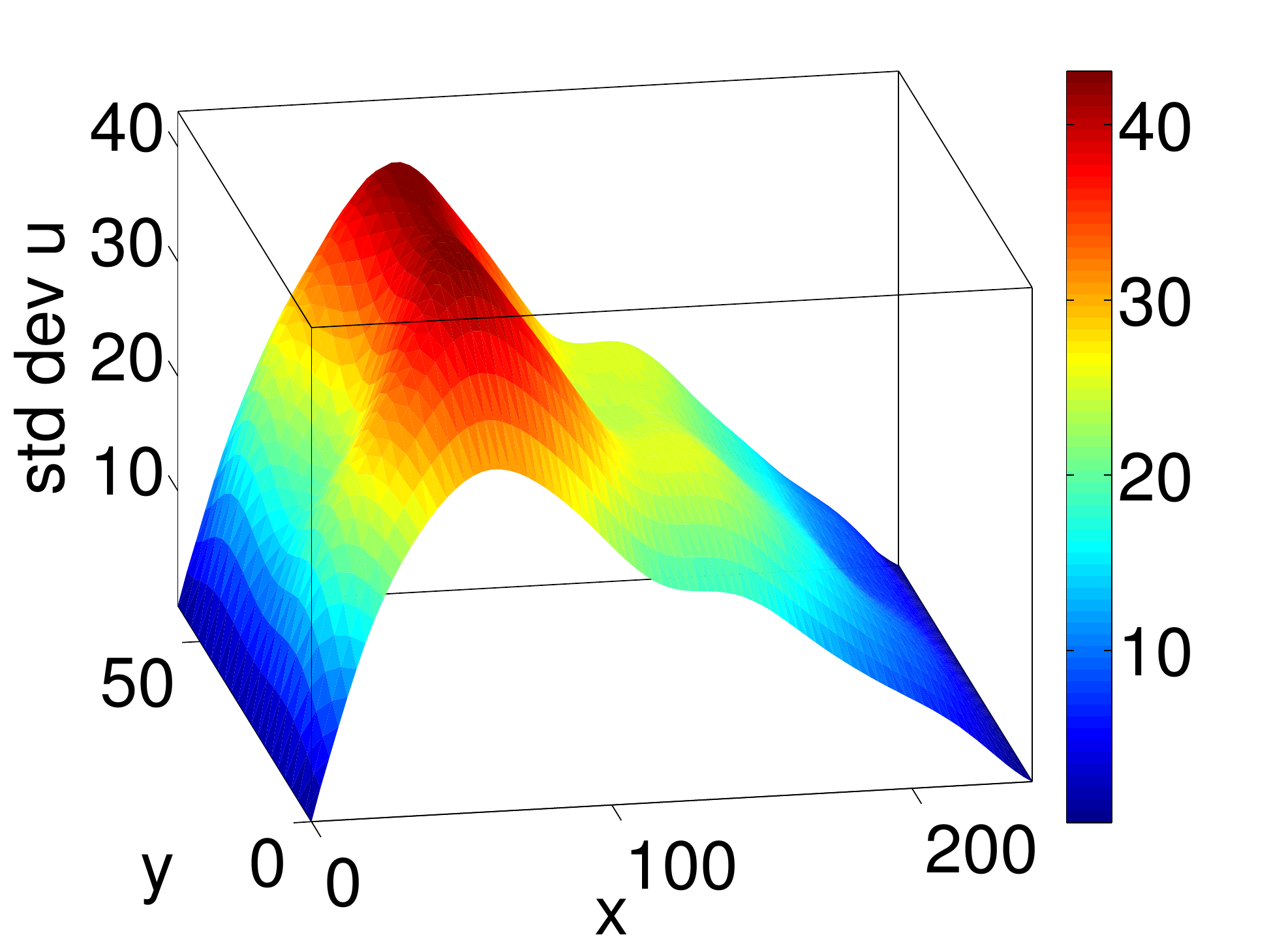}
        \caption{Std. Dev, $\xi, d=10$}
    \end{subfigure}        
    \begin{subfigure}[t]{0.3\textwidth}
        \centering
        \includegraphics[height=1.2in]{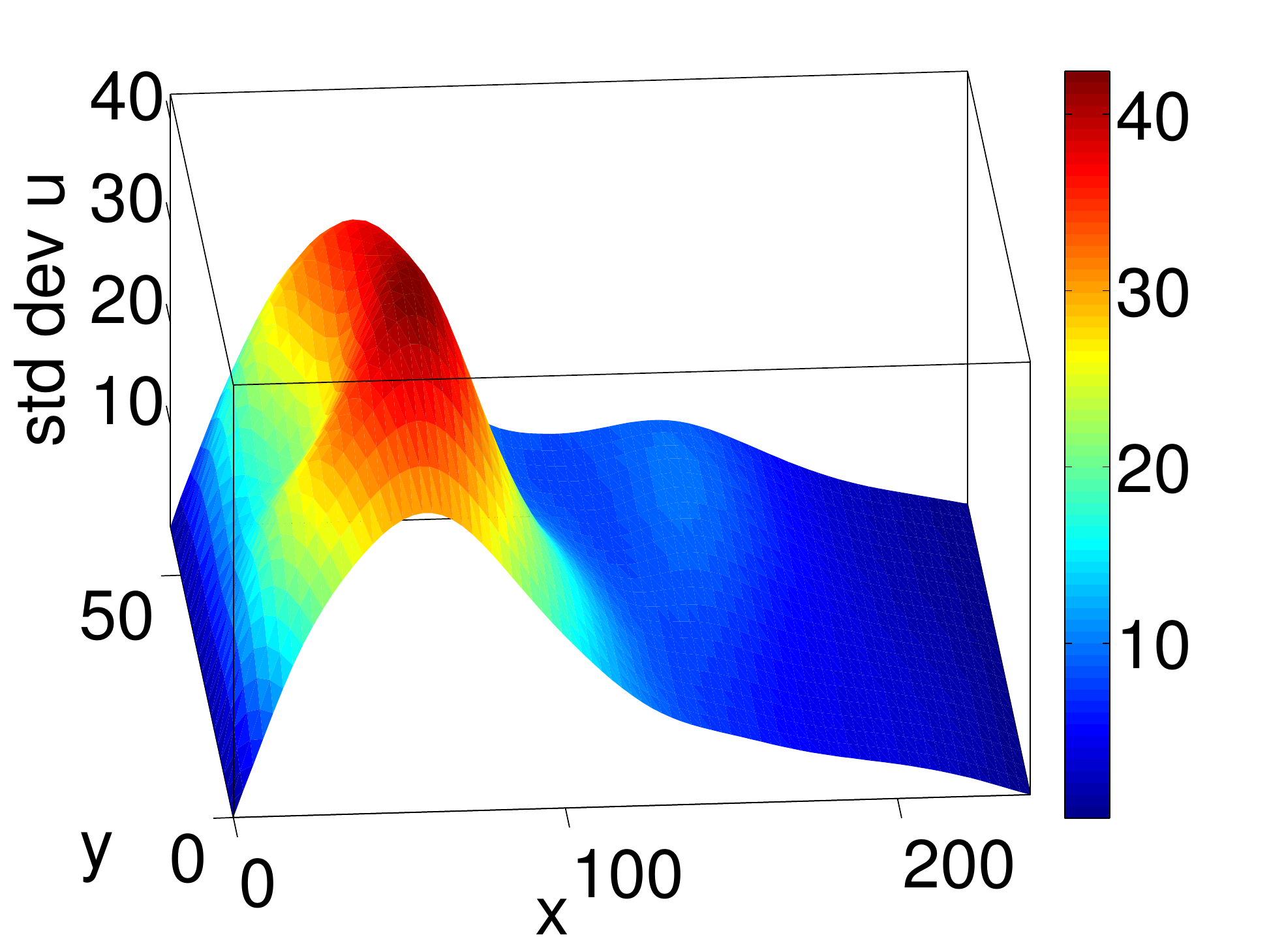}
        \caption{Std. Dev, $\eta, d=3$}
    \end{subfigure}    
    \begin{subfigure}[t]{0.3\textwidth}
        \centering
        \includegraphics[height=1.2in]{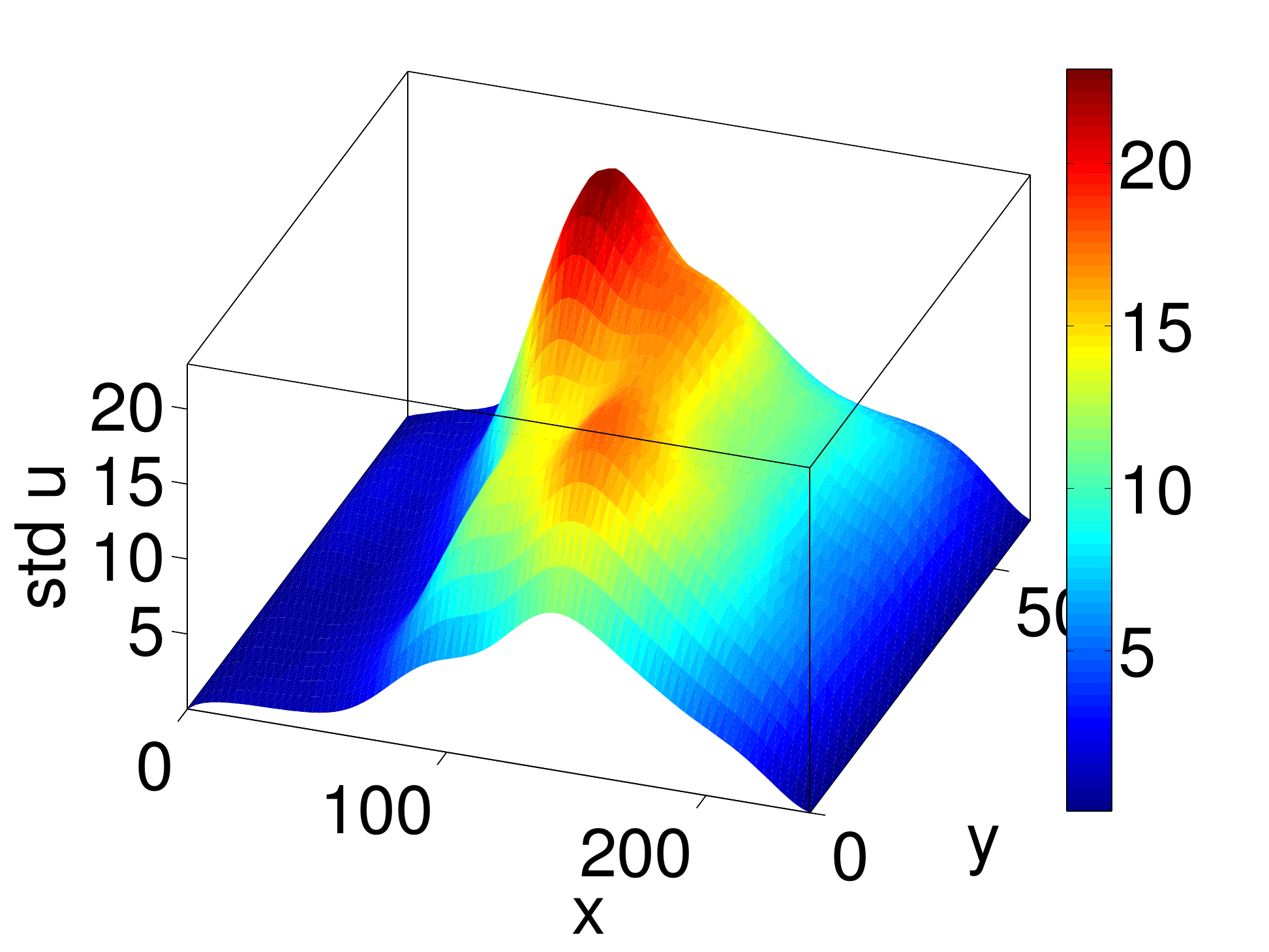}
        \caption{Error, $\eta, d=3$}
    \end{subfigure}  
    
    \caption{Standard deviation of the solution with the stochastic basis adapted for subdomain $D_1$ with random variables, $\eta$ and dimension = 3, order = 3.} \label{RK:fig:u_std_D1}
\end{figure}

\begin{figure}[t!]
    \centering
    \begin{subfigure}[t]{0.3\textwidth}
        \centering
       \includegraphics[height=1.2in]{figures/u_std_xid10_p3}
        \caption{Std. Dev, $\xi, d=10$}
    \end{subfigure}        
    \begin{subfigure}[t]{0.3\textwidth}
        \centering
        \includegraphics[height=1.2in]{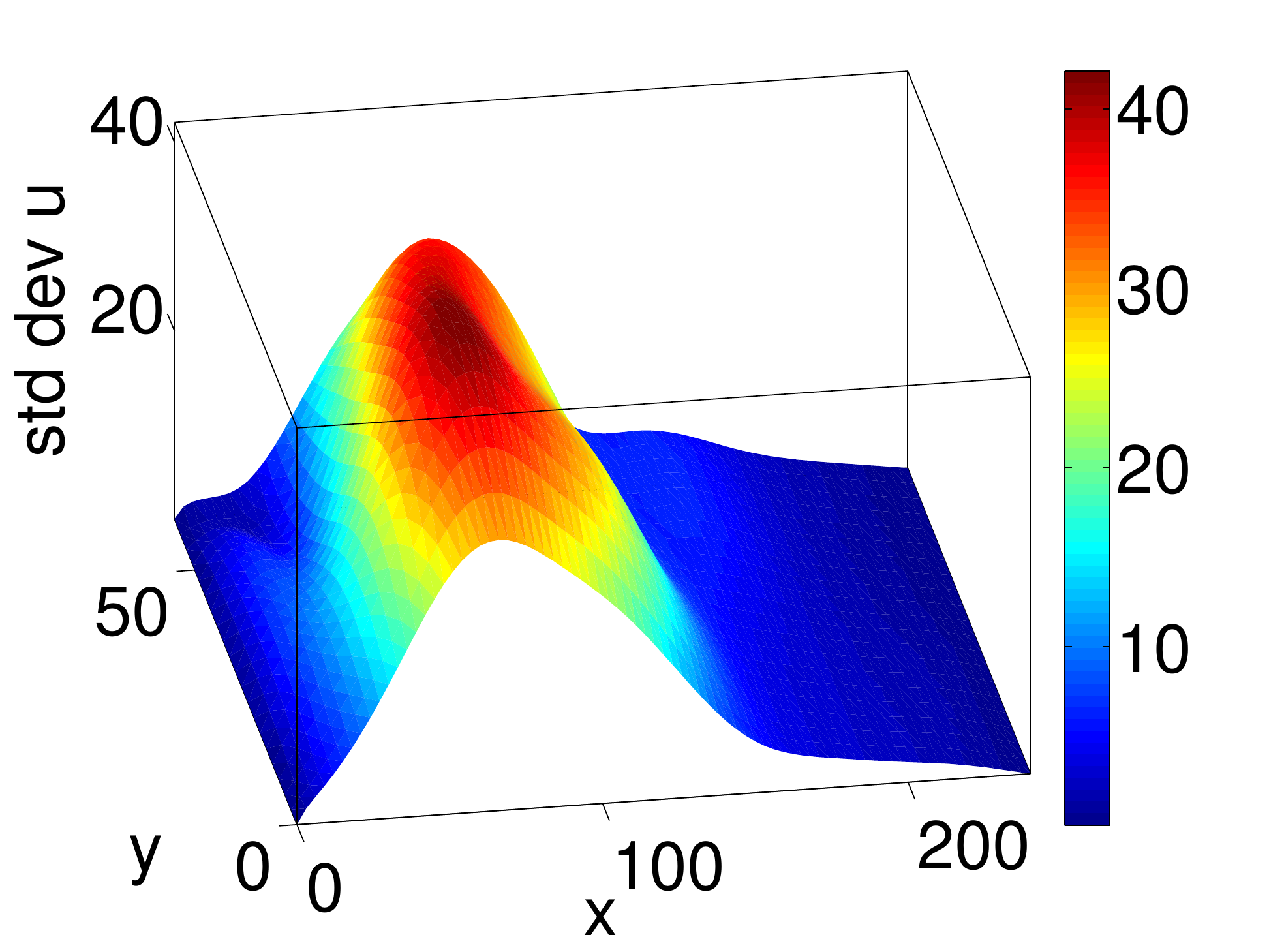}
        \caption{Std. Dev, $\eta, d=3$}
    \end{subfigure}    
    \begin{subfigure}[t]{0.3\textwidth}
        \centering
        \includegraphics[height=1.2in]{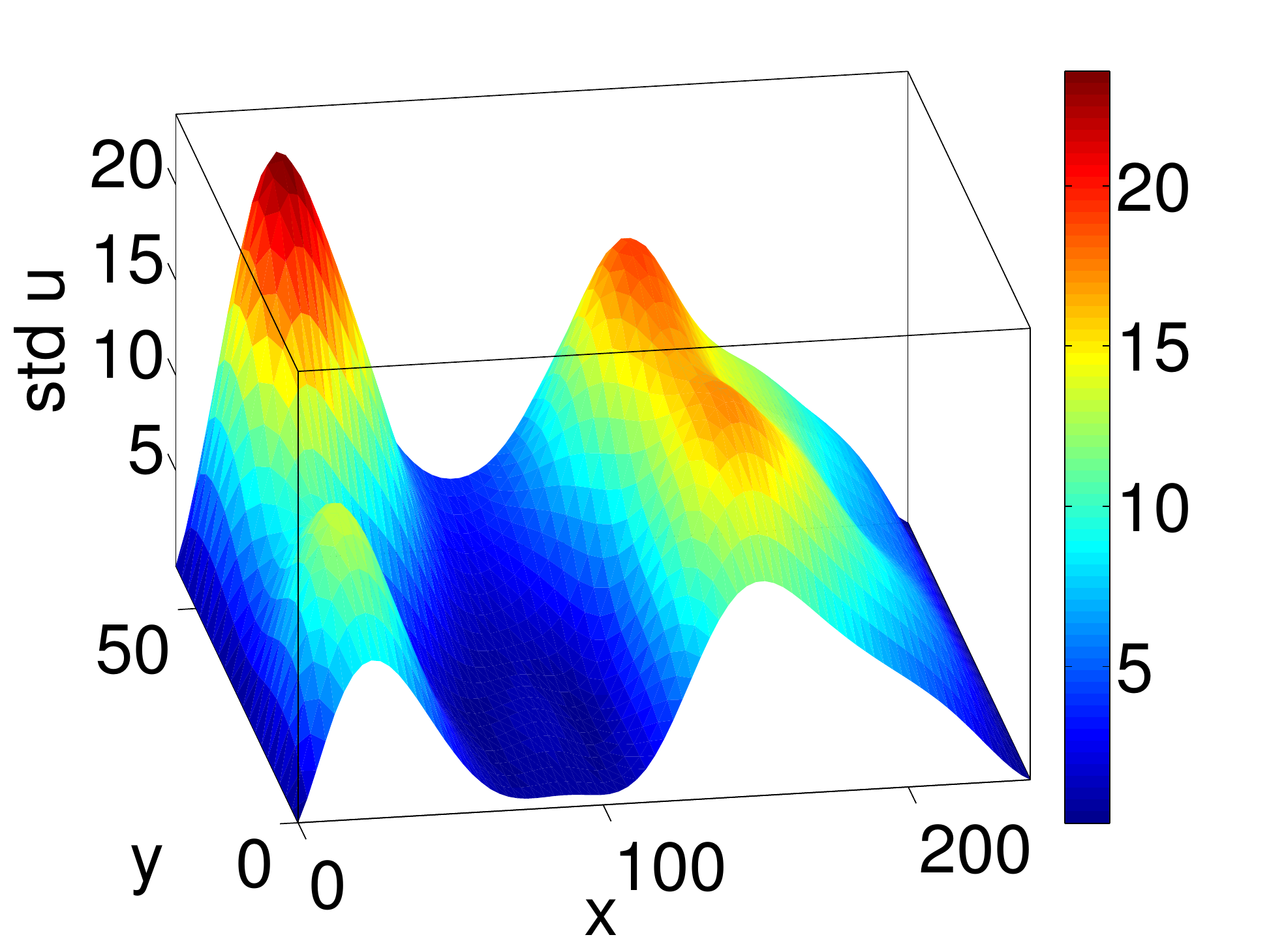}
        \caption{Error, $\eta, d=3$}
    \end{subfigure}  
    
    \caption{Standard deviation of the solution with the stochastic basis adapted for subdomain $D_2$ with random variables, $\eta$ and dimension = 3, order = 3.} \label{RK:fig:u_std_D2}
\end{figure}

Figures \ref{RK:fig:u_mean_D1}-\ref{RK:fig:u_std_D2} present the mean and standard deviation in subdomains $D_1$ and $D_2$, computed with full dimension ($d=10$) and reduced dimension ($d$=3) for the first type of boundary conditions. For example, Figure \ref{RK:fig:u_mean_D1_full} corresponds to the mean of the solution computed with full dimension $d=10$ and sparse-grid level 5. Figure~\ref{RK:fig:u_mean_D1_reduced} shows the mean of the solution computed with the low-dimensional  $d=3$ basis adapted to the subdomain $D_1.$ Similar behavior is observed in other subdomains. 

\begin{figure}[t!]
    \centering
    \begin{subfigure}[t]{0.3\textwidth}
        \centering
        \includegraphics[height=1.2in]{figures/u_mean_xid10_p3}
        \caption{Mean, $\xi, d=10$}
    \end{subfigure}        
    \begin{subfigure}[t]{0.3\textwidth}
        \centering
        \includegraphics[height=1.2in]{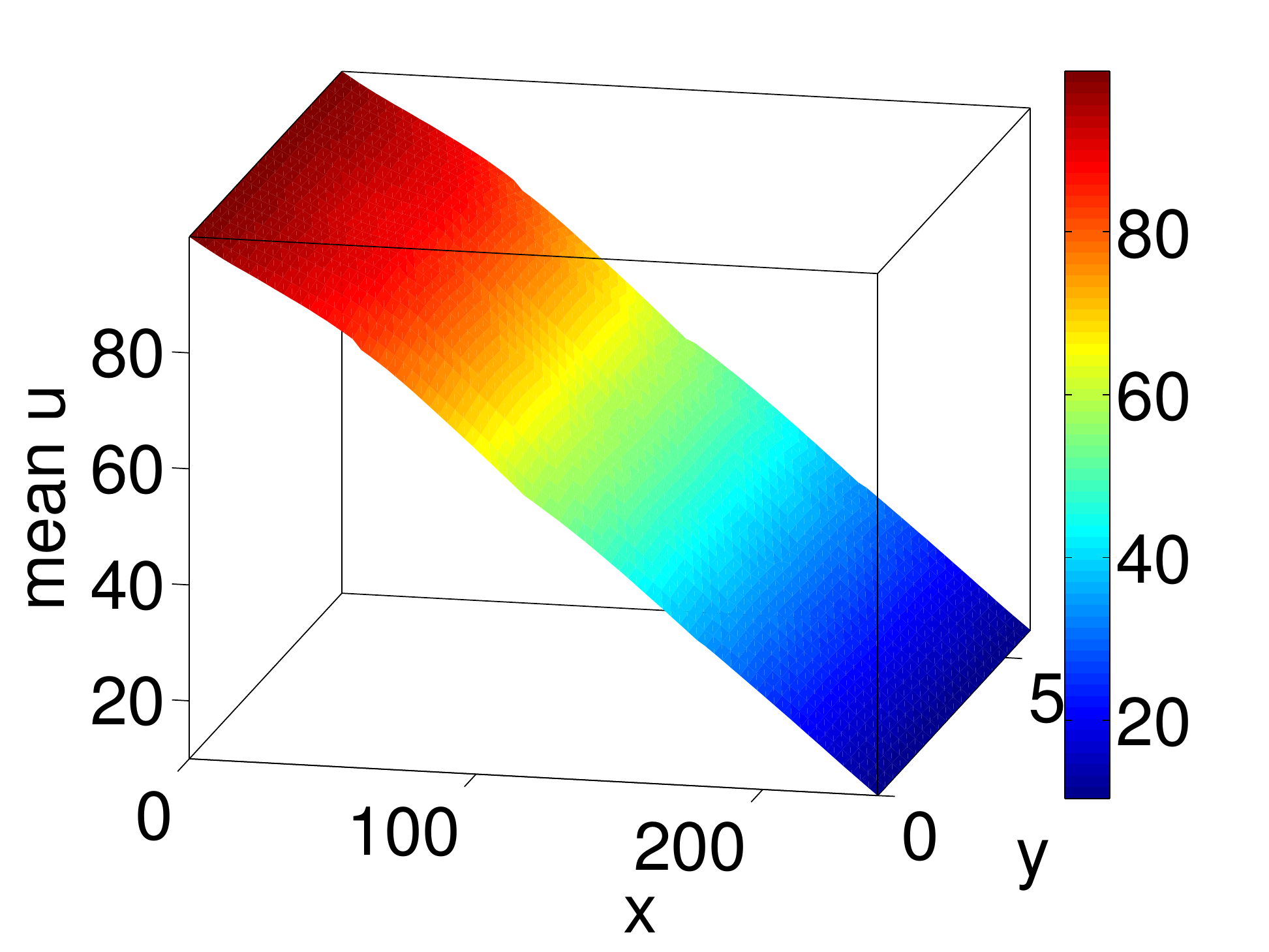}
        \caption{Mean, $\eta, d=3$}
    \end{subfigure}    
    \begin{subfigure}[t]{0.3\textwidth}
        \centering
        \includegraphics[height=1.2in]{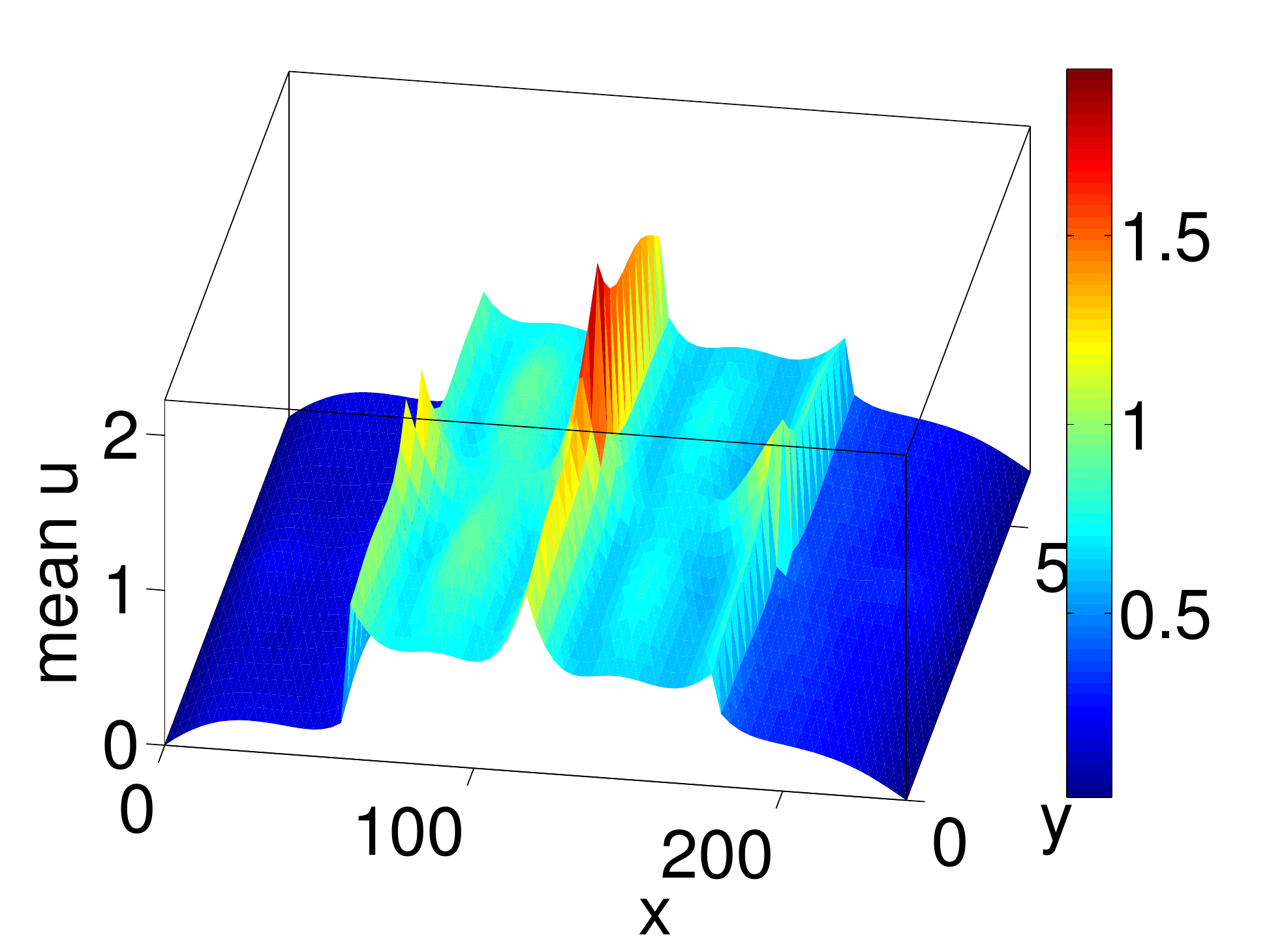}
        \caption{Error, $\eta, d=3$}
    \end{subfigure}  

    \caption{Mean of the solution in entire domain $D$ for \eqref{RK:eq:spde}.} \label{RK:fig:u_mean_combine}
\end{figure}

\begin{figure}[t!]
    \centering
    \begin{subfigure}[t]{0.3\textwidth}
        \centering
        \includegraphics[height=1.2in]{figures/u_std_xid10_p3}
        \caption{Std. Dev, $\xi, d=10$}
    \end{subfigure}        
    \begin{subfigure}[t]{0.3\textwidth}
        \centering
        \includegraphics[height=1.2in]{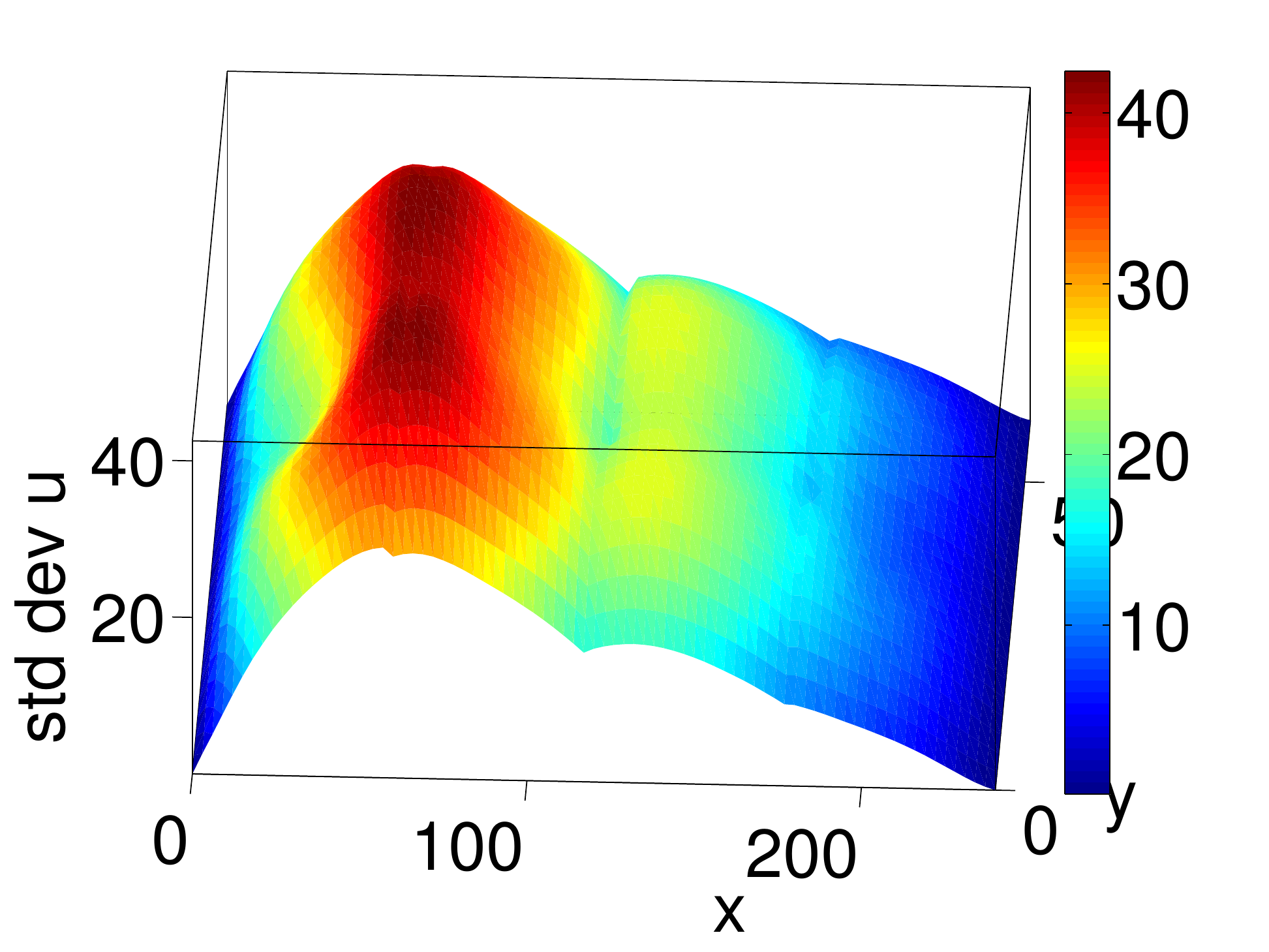}
        \caption{Std. Dev, $\eta, d=3$}
    \end{subfigure}    
    \begin{subfigure}[t]{0.3\textwidth}
        \centering
        \includegraphics[height=1.2in]{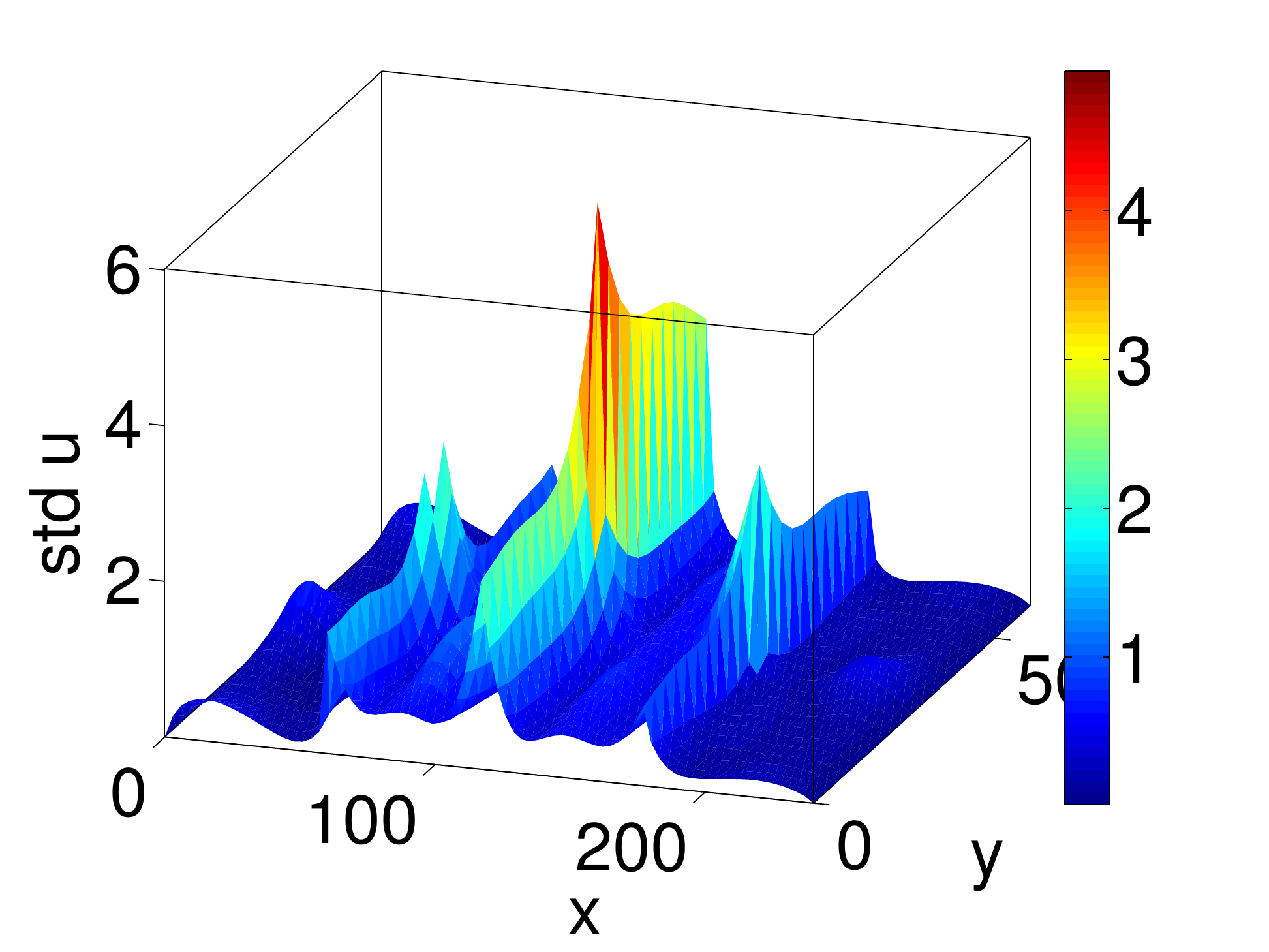}
        \caption{Error, $\eta, d=3$}
    \end{subfigure}  
    
    \caption{Standard deviation of the solution in entire domain $D$ for \eqref{RK:eq:spde}.} \label{RK:fig:u_std_combine}
\end{figure}

Figure~\ref{RK:fig:u_mean_D1_error} shows the error in the mean. Note that although the basis is adapted only to subdomain $D_1$, the low-dimensional solution is computed in the entire spatial domain $D$.  Hence, in Figure~\ref{RK:fig:u_mean_D1}, we are only interested in the solution in subdomain $D_1.$ We can see the error in subdomain $D_1$ is quite small compared to other domains. Similarly, Figure~\ref{RK:fig:u_mean_D2} depicts a comparison of the mean solution in subdomain $D_2$. A similar comparison for the standard deviation is shown in Figures~\ref{RK:fig:u_std_D1} and \ref{RK:fig:u_std_D2}. We also computed and collected the adapted solution from each subdomain and plotted the mean and standard deviations in Figures~\ref{RK:fig:u_mean_combine} and~\ref{RK:fig:u_std_combine} respectively. 

\begin{figure}[t!]
    \centering
    \begin{subfigure}[t]{0.3\textwidth}
        \centering
        \includegraphics[height=1.2in]{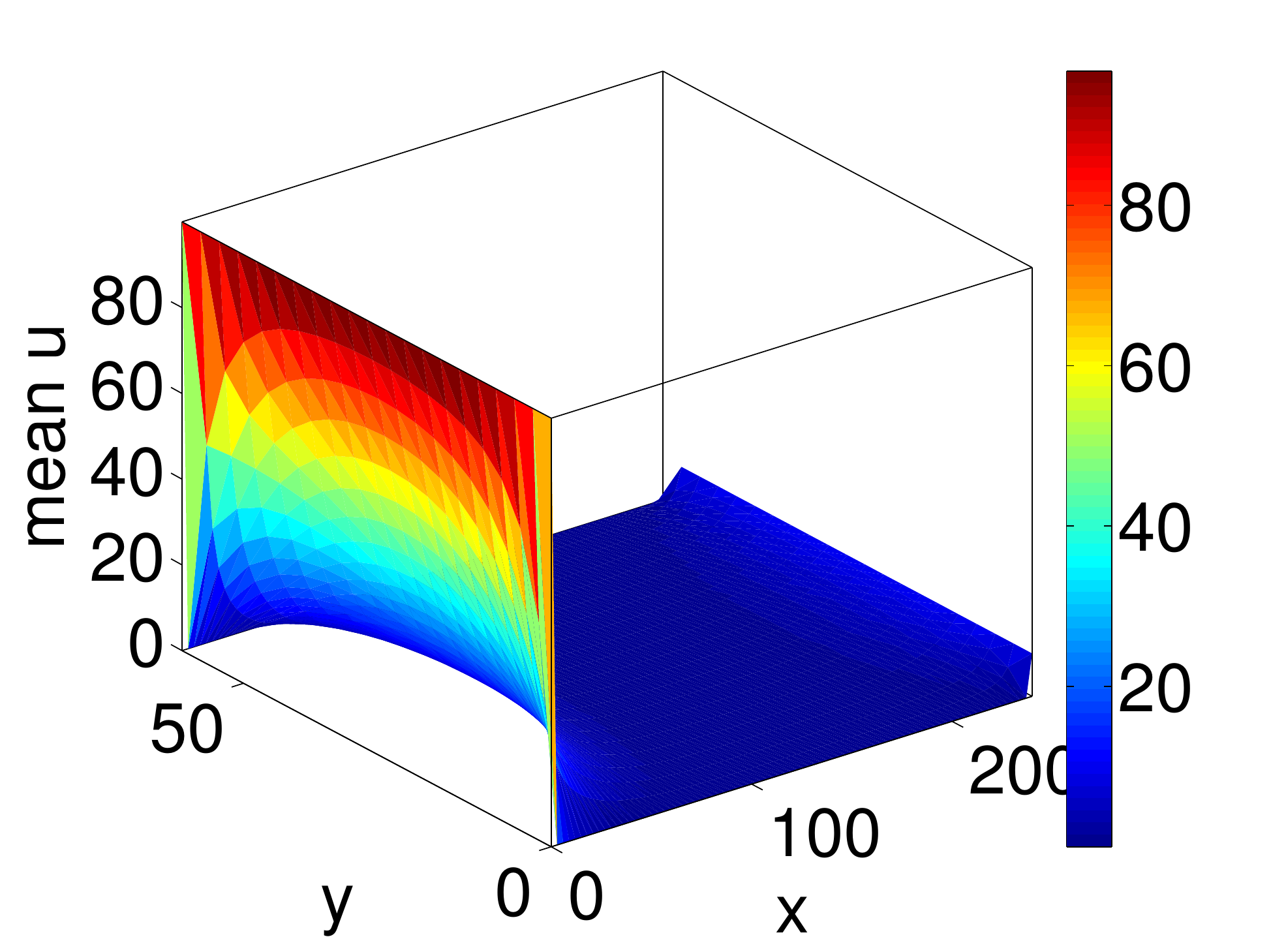}
        \caption{Mean, $\xi, d=10$}
    \end{subfigure}        
    \begin{subfigure}[t]{0.3\textwidth}
        \centering
        \includegraphics[height=1.2in]{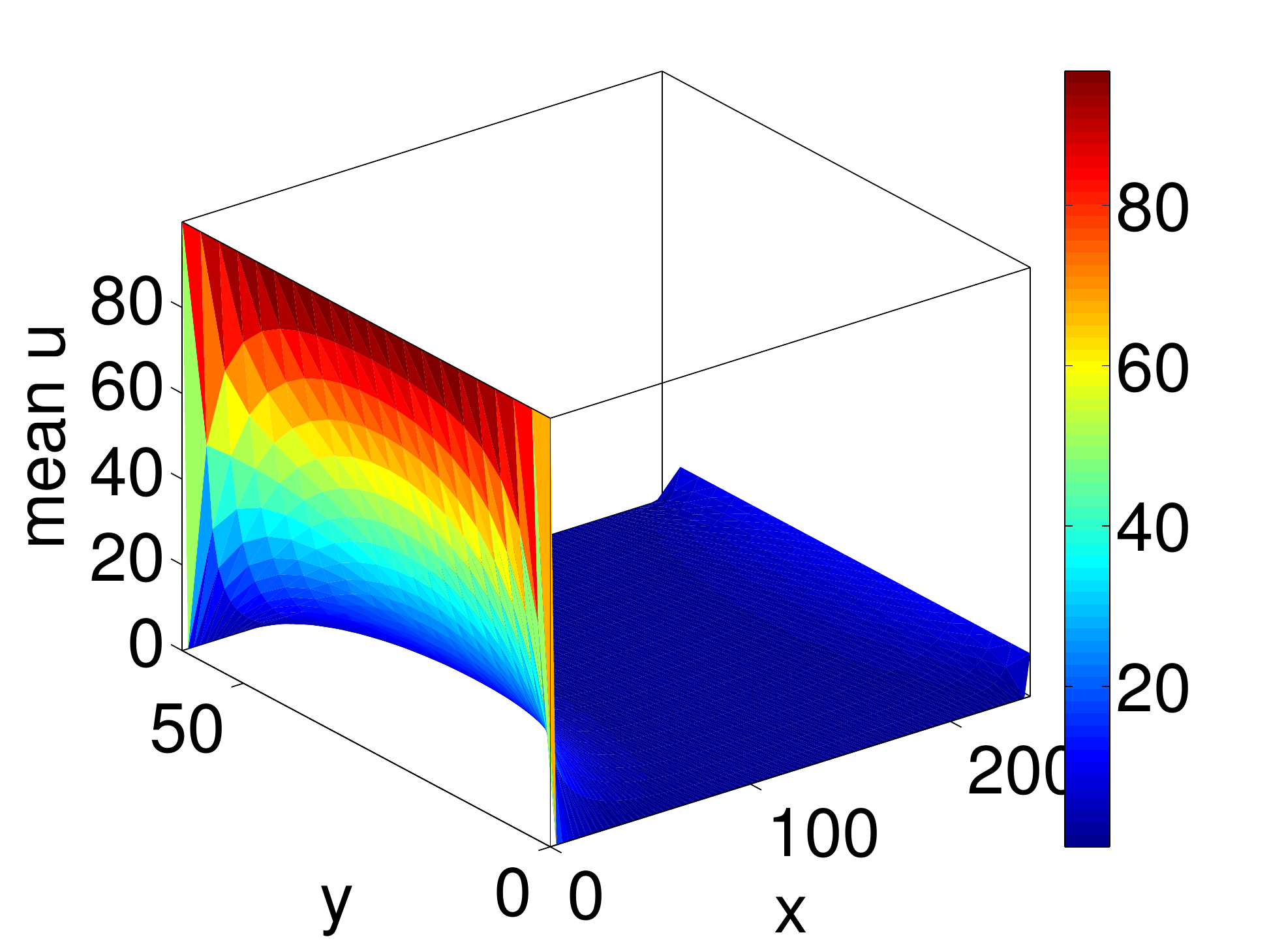}
        \caption{Mean, $\eta, d=3$}
    \end{subfigure}    
    \begin{subfigure}[t]{0.3\textwidth}
        \centering
        \includegraphics[height=1.2in]{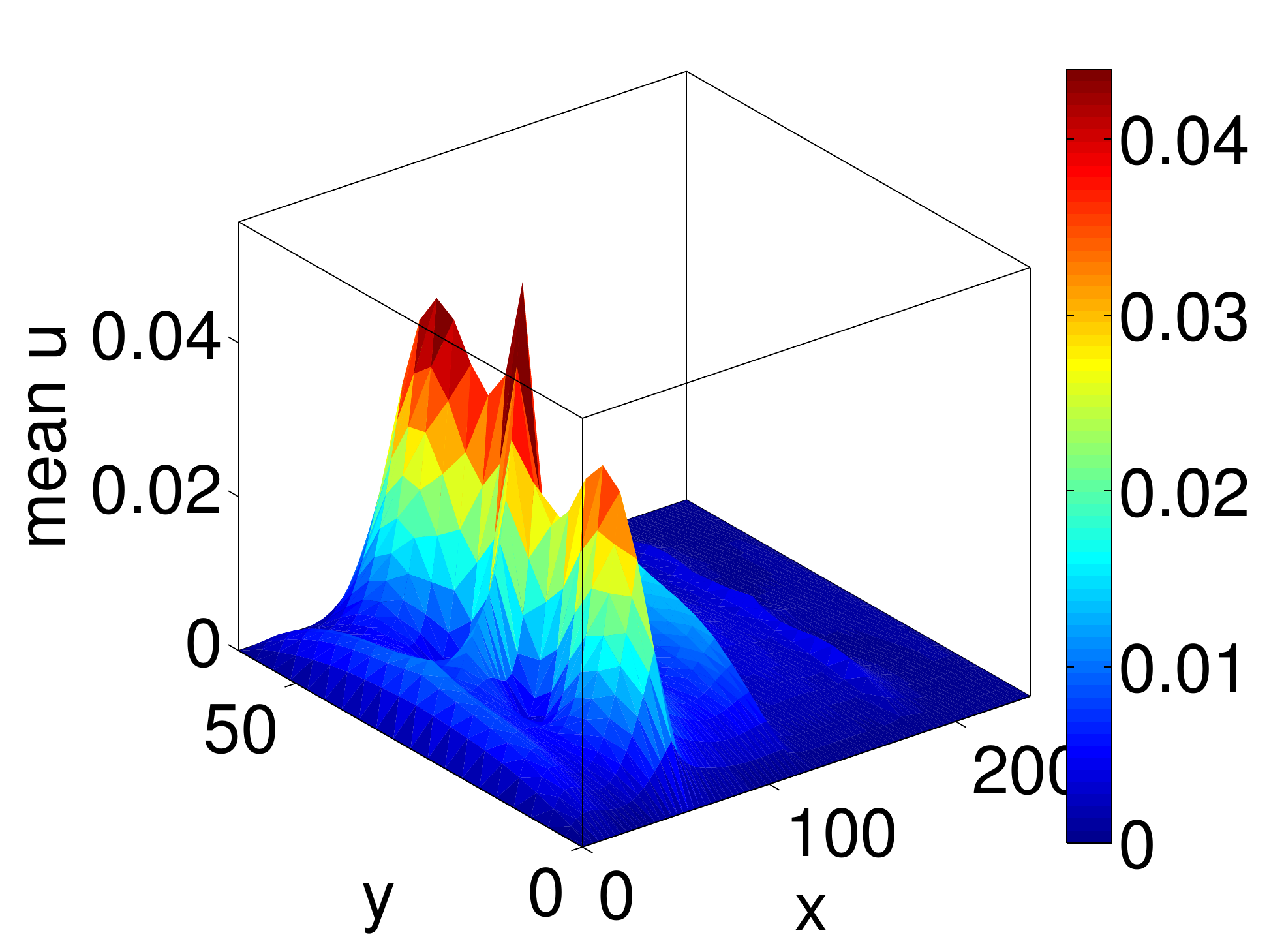}
        \caption{Error, $\eta, d=3$}
    \end{subfigure}  

    \caption{Mean of the solution in entire domain $D$ for \eqref{RK:eq:spde_dbc}.} \label{RK:fig:u_mean_combine_DBC}
\end{figure}

\begin{figure}[t!]
    \centering
    \begin{subfigure}[t]{0.3\textwidth}
        \centering
        \includegraphics[height=1.2in]{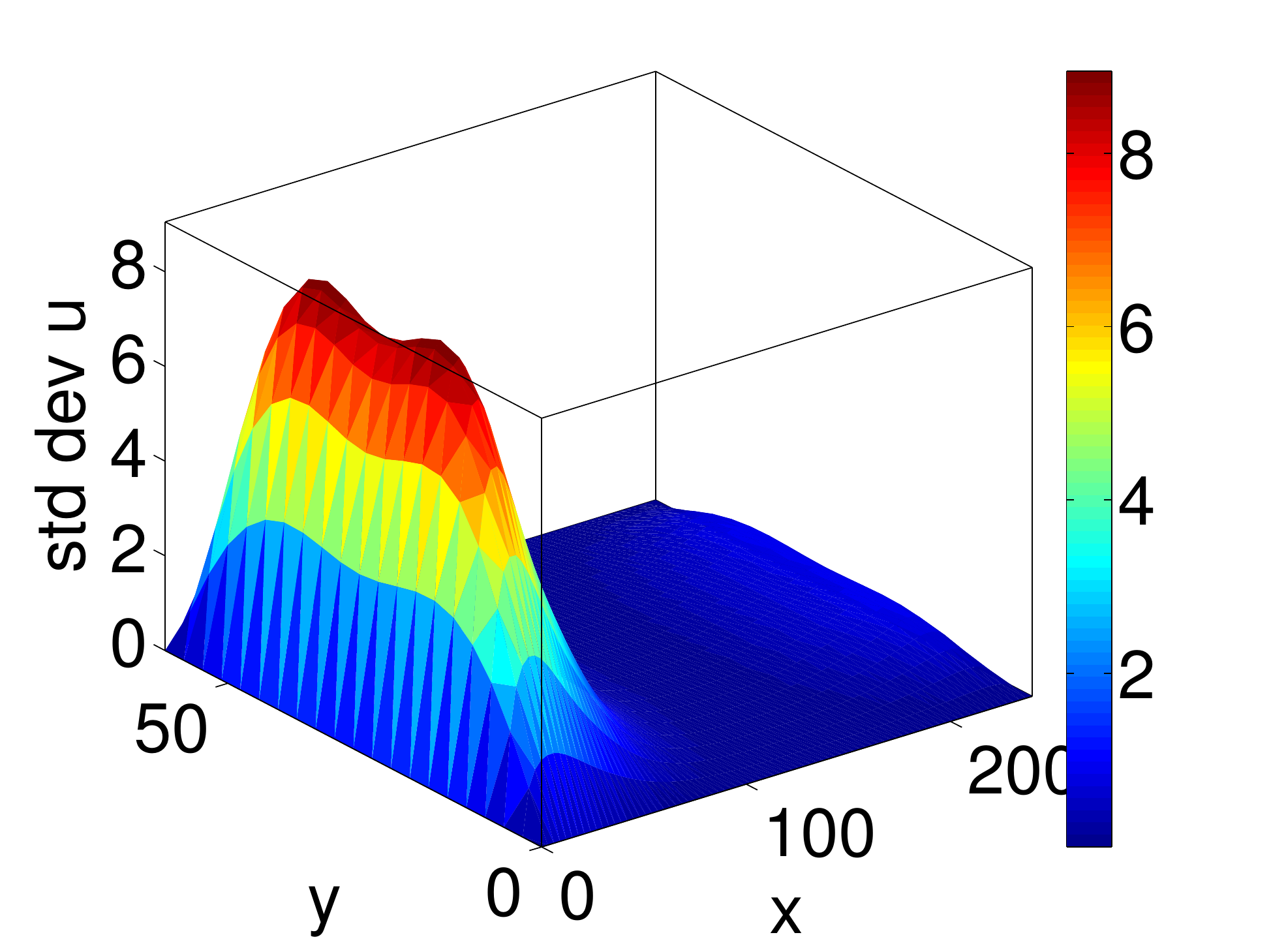}
        \caption{Std. Dev, $\xi, d=10$}
    \end{subfigure}        
    \begin{subfigure}[t]{0.3\textwidth}
        \centering
       \includegraphics[height=1.2in]{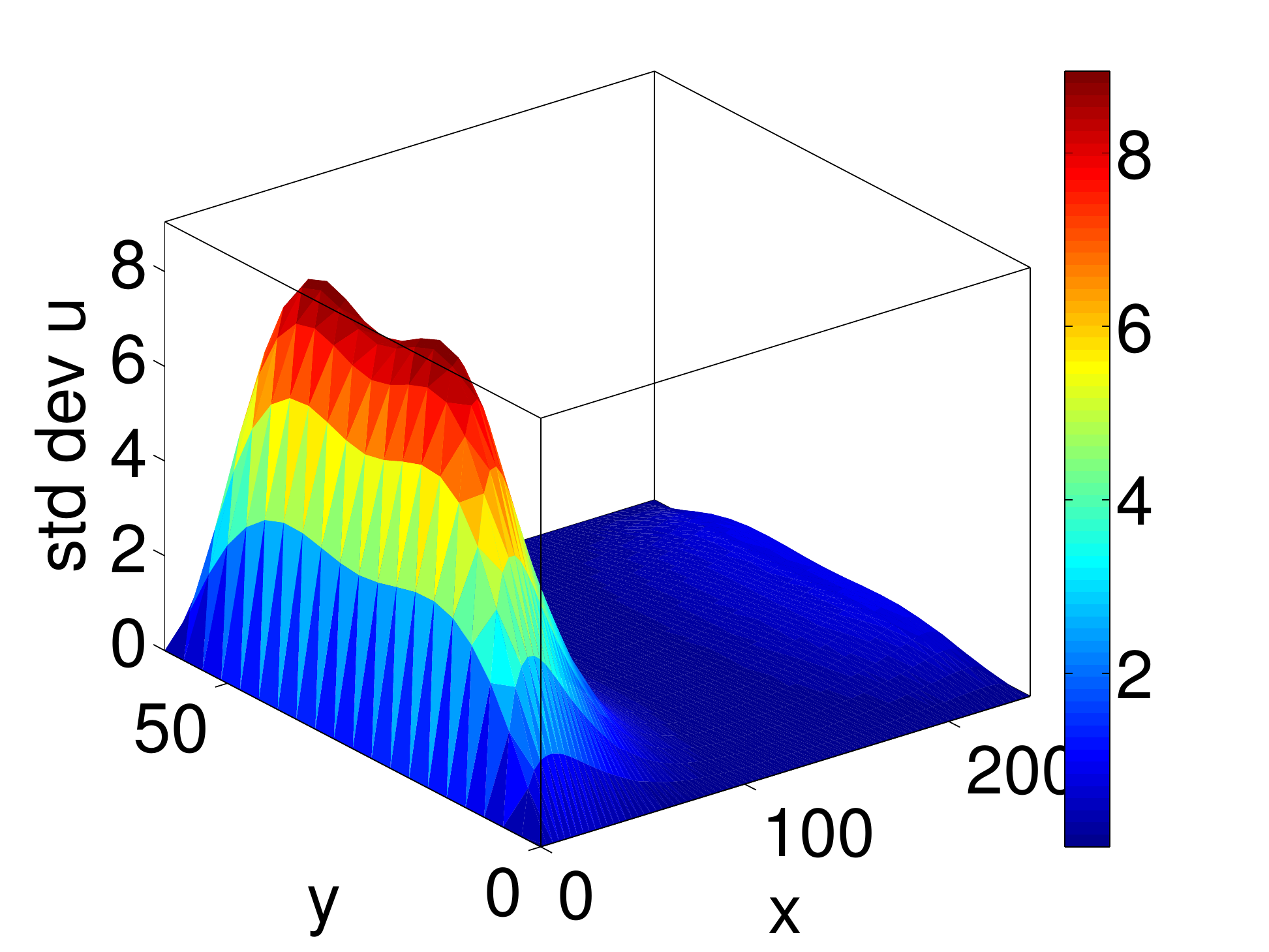}
        \caption{Std. Dev, $\eta, d=3$}
    \end{subfigure}    
    \begin{subfigure}[t]{0.3\textwidth}
        \centering
        \includegraphics[height=1.2in]{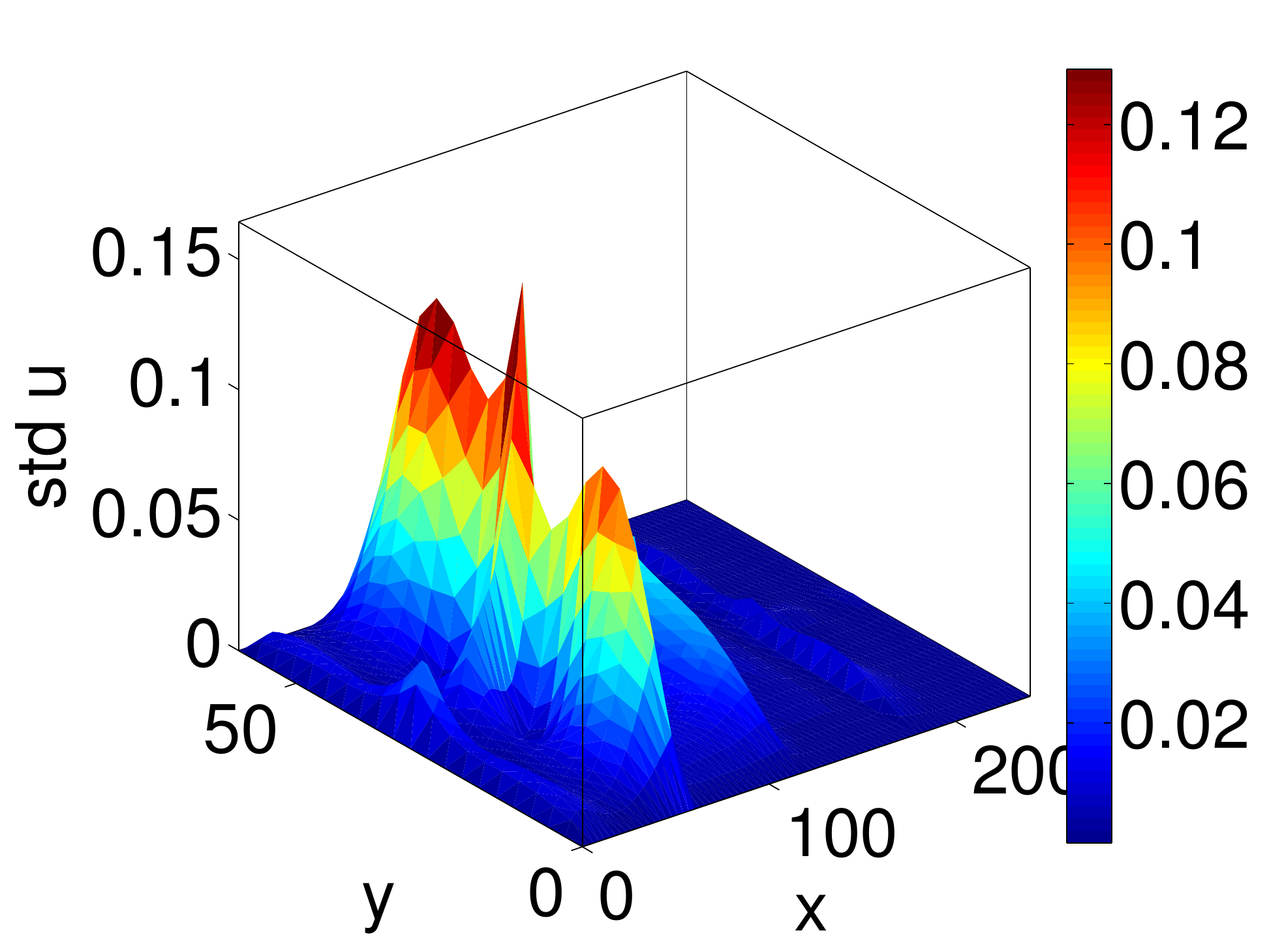}
        \caption{Error, $\eta, d=3$}
    \end{subfigure}  
    
    \caption{Standard deviation of the solution in entire domain $D$ for \eqref{RK:eq:spde_dbc}.} \label{RK:fig:u_std_combine_DBC}
\end{figure}

For the second case (Dirichlet boundary conditions on all sides), we show the mean and standard deviation in Figures \ref{RK:fig:u_mean_combine_DBC} and \ref{RK:fig:u_std_combine_DBC}. We observe very good agreement for the mean and standard deviation results obtained with full and reduced dimensions. We also recognize a small mismatch of the solutions from adjacent subdomains at the subdomain boundaries. This discrepancy is due to the low-dimensional approximation, which uses a different adapted basis in each subdomain. If required, this inconsistency can be eliminated by appropriate post-processing, such as averaging, interpolation, or projecting each low-dimensional solution onto a global basis. 

\begin{figure}[t!]
    \centering
      \begin{subfigure}[t]{0.3\textwidth}
        \centering
        \includegraphics[height=1.2in]{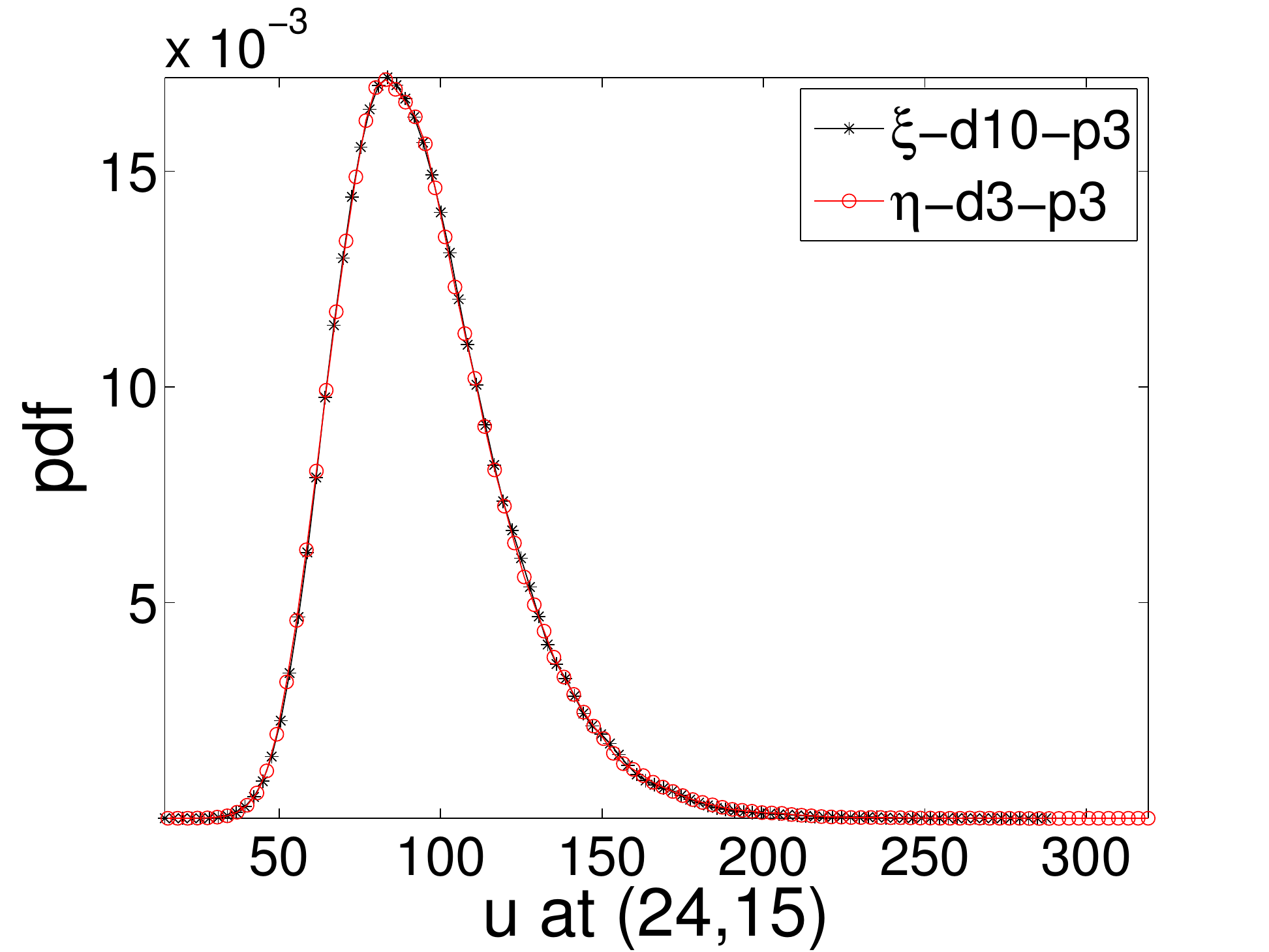}
        \caption{pdf in $D_1$, $\eta, d=3$}
     \end{subfigure}
    \begin{subfigure}[t]{0.3\textwidth}
        \centering
        \includegraphics[height=1.2in]{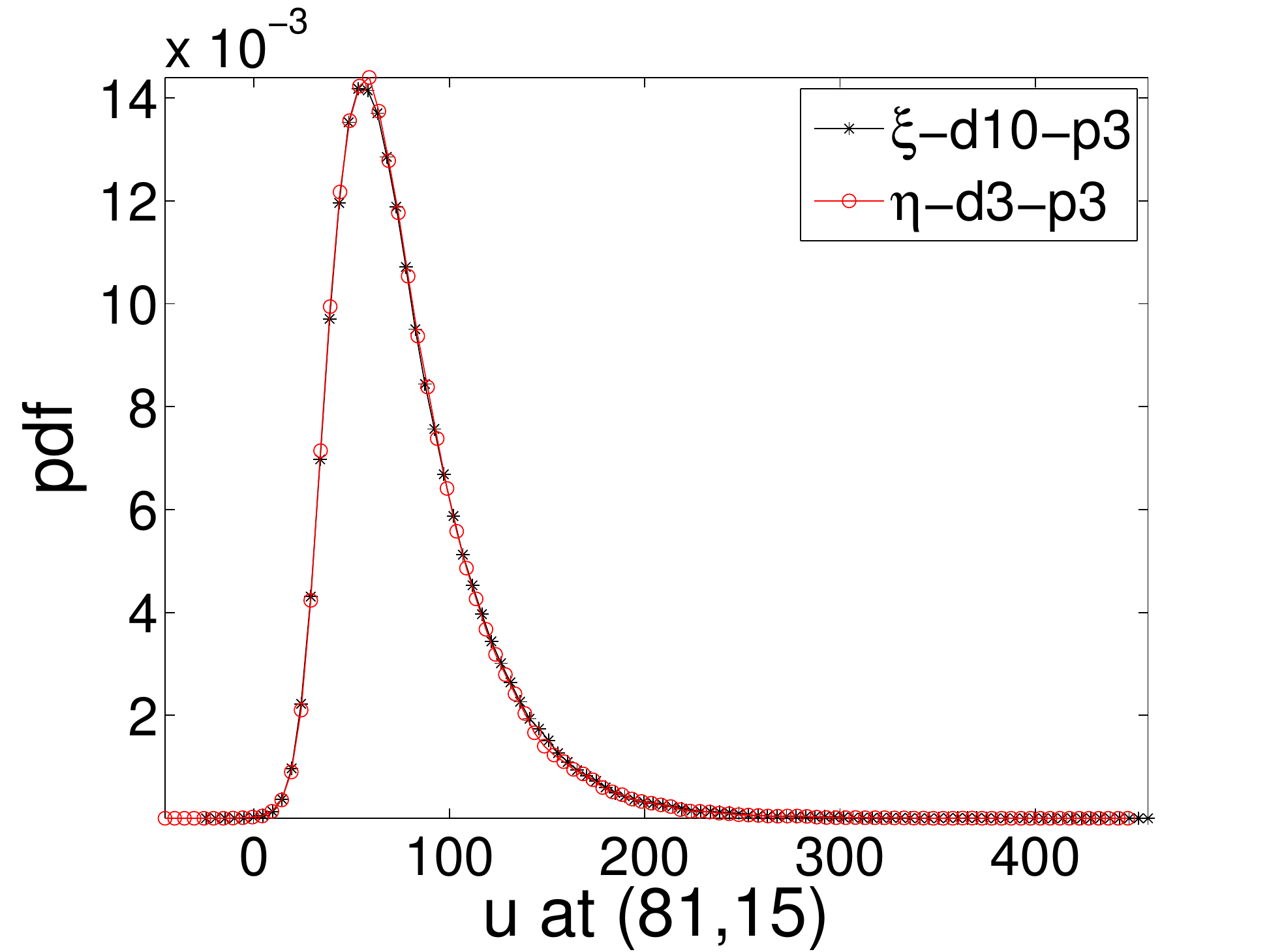}
        \caption{pdf in $D_2$, $\eta, d=3$}
     \end{subfigure}
     \begin{subfigure}[t]{0.3\textwidth}
        \centering
        \includegraphics[height=1.2in]{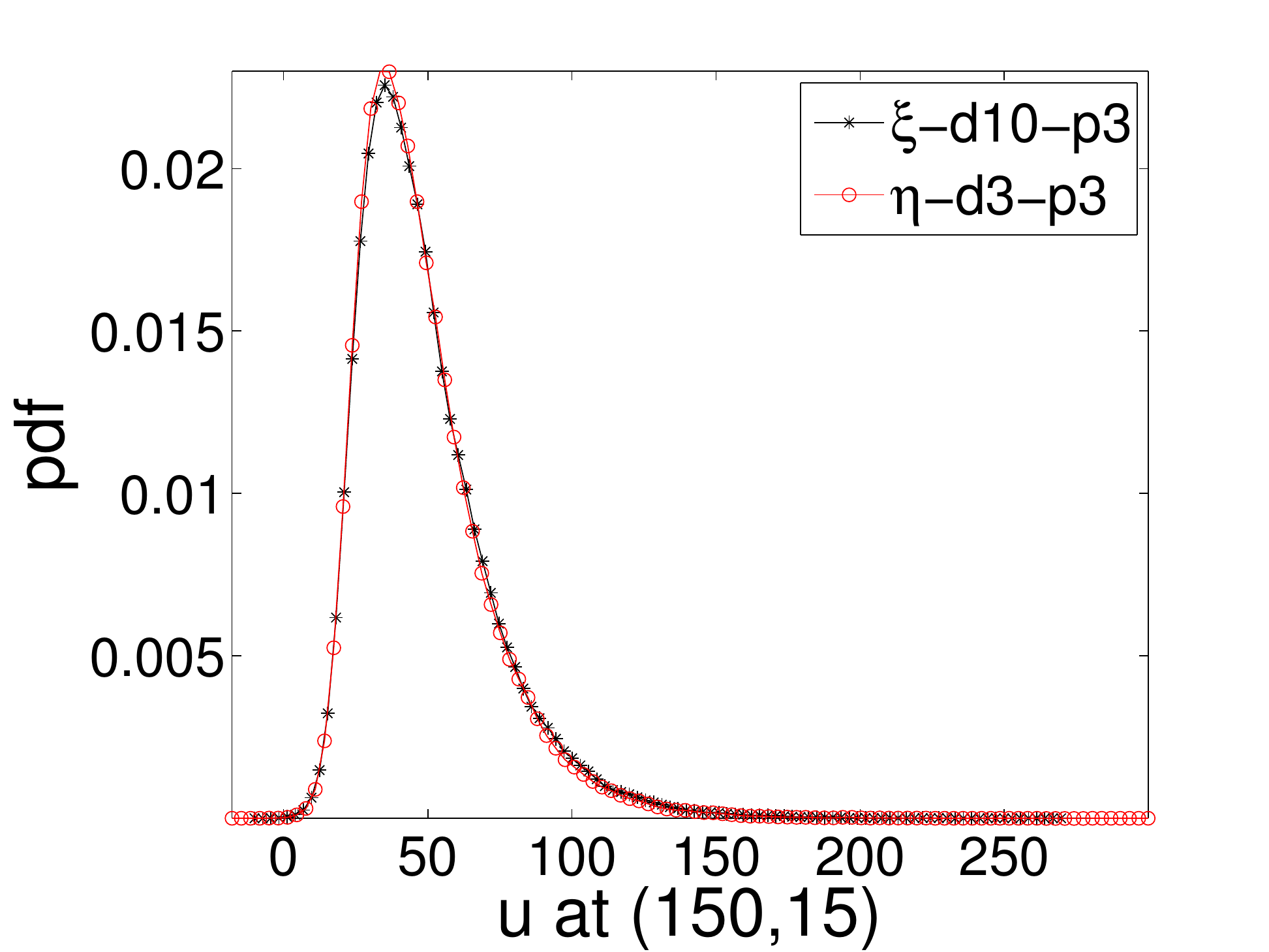}
        \caption{pdf in $D_3$, $\eta, d=3$}
     \end{subfigure}
    \caption{Probability density function of the solution at points (24,15) in $D_1$, (81, 15) in $D_2$ and (150,15) in $D_3.$} \label{RK:fig:pdf_D123}
\end{figure}

\begin{figure}[t!]
    \centering
      \begin{subfigure}[t]{0.3\textwidth}
        \centering
        \includegraphics[height=1.2in]{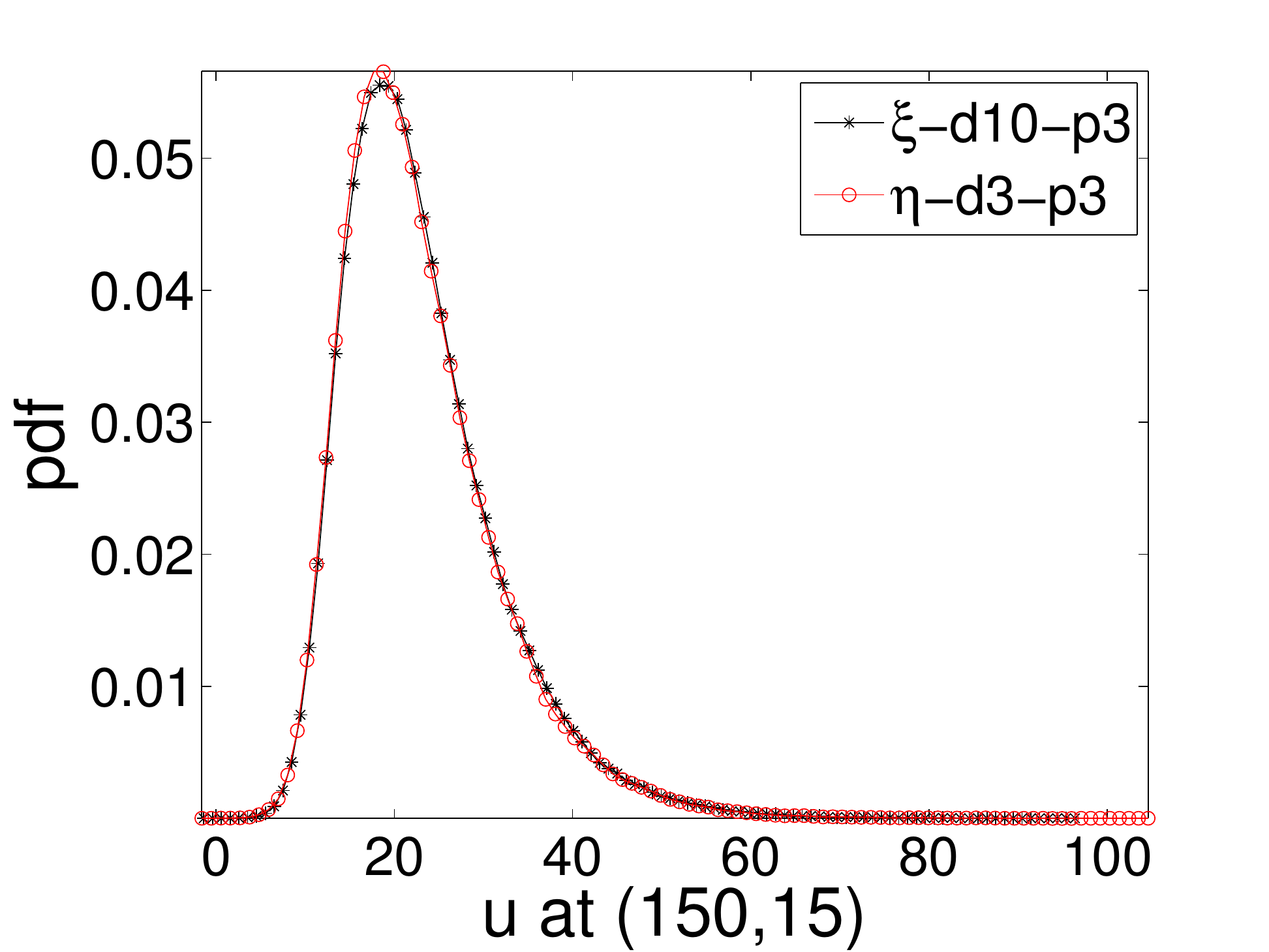}
        \caption{pdf in $D_4$, $\eta, d=3$}
     \end{subfigure}
    \begin{subfigure}[t]{0.3\textwidth}
        \centering
        \includegraphics[height=1.2in]{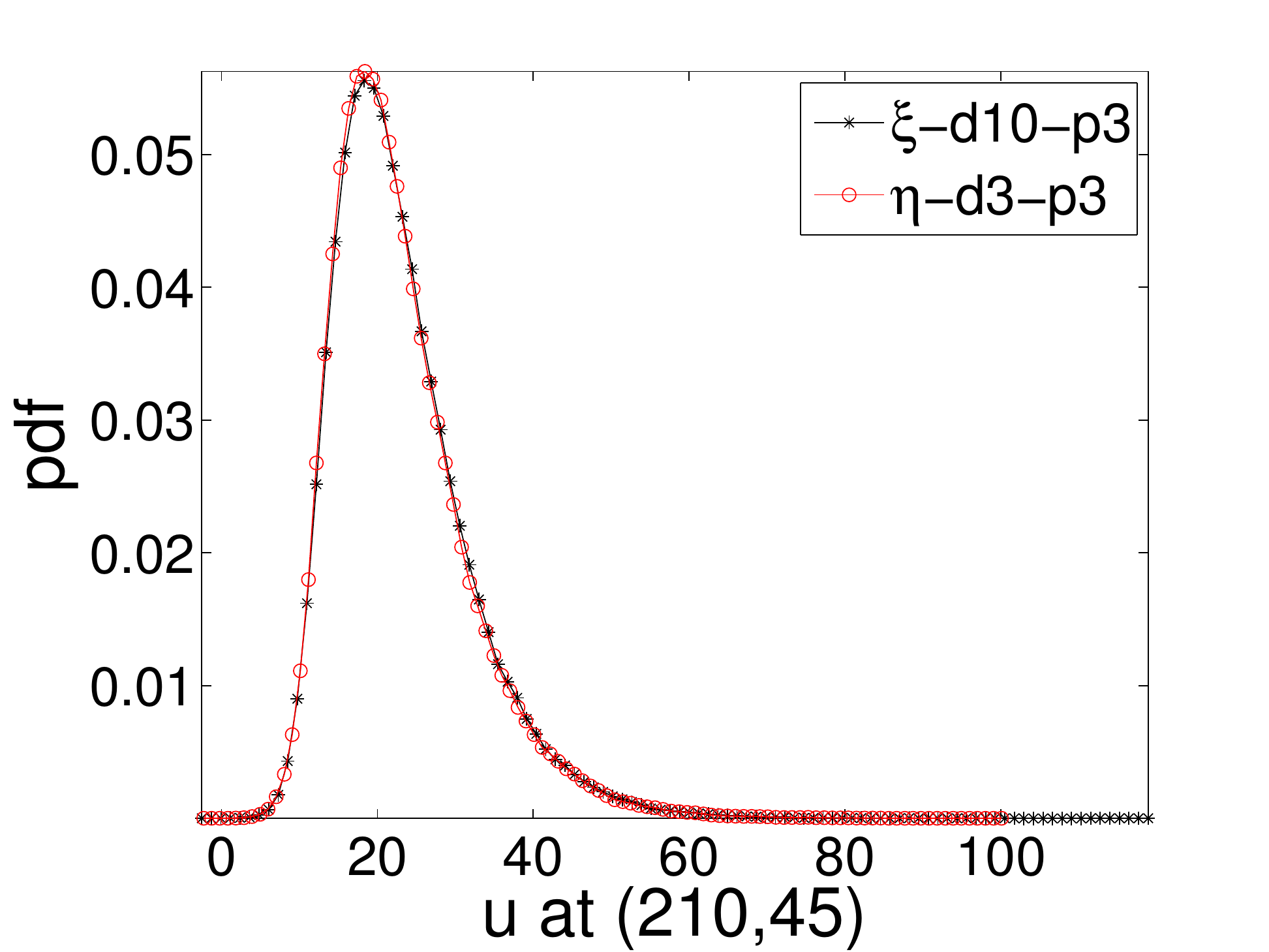}
        \caption{pdf in $D_5$, $\eta, d=3$}
     \end{subfigure}
     \begin{subfigure}[t]{0.3\textwidth}
        \centering
        \includegraphics[height=1.2in]{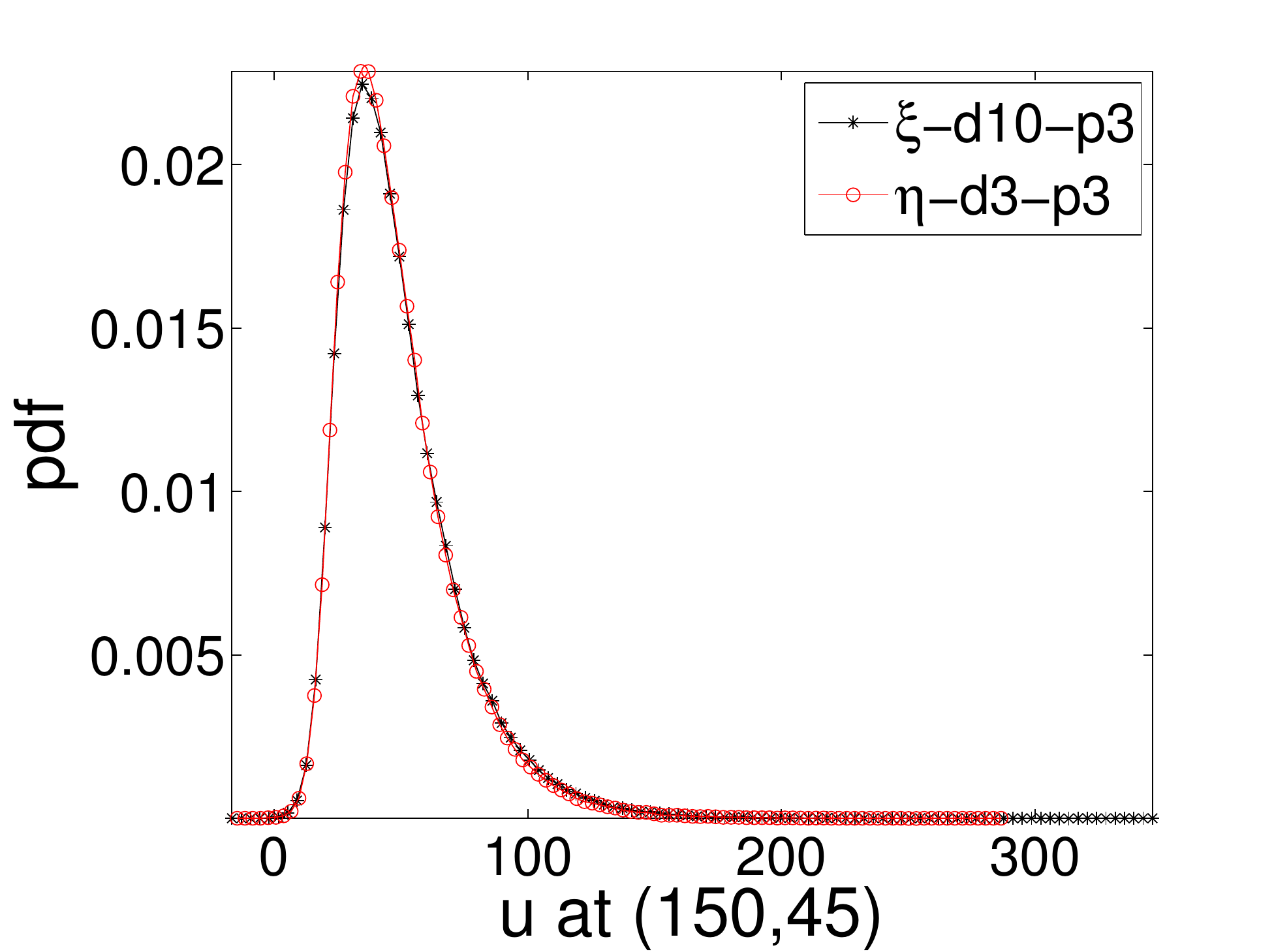}
        \caption{pdf in $D_6$, $\eta, d=3$}
     \end{subfigure}
    \caption{Probability density function of the solution at points (210,15) in $D_4$, (210, 45) in $D_5$ and (150,45) in $D_6.$} \label{RK:fig:pdf_D456}
\end{figure}

\begin{figure}[t!]
    \centering
      \begin{subfigure}[t]{0.45\textwidth}
        \centering
        \includegraphics[height=1.2in]{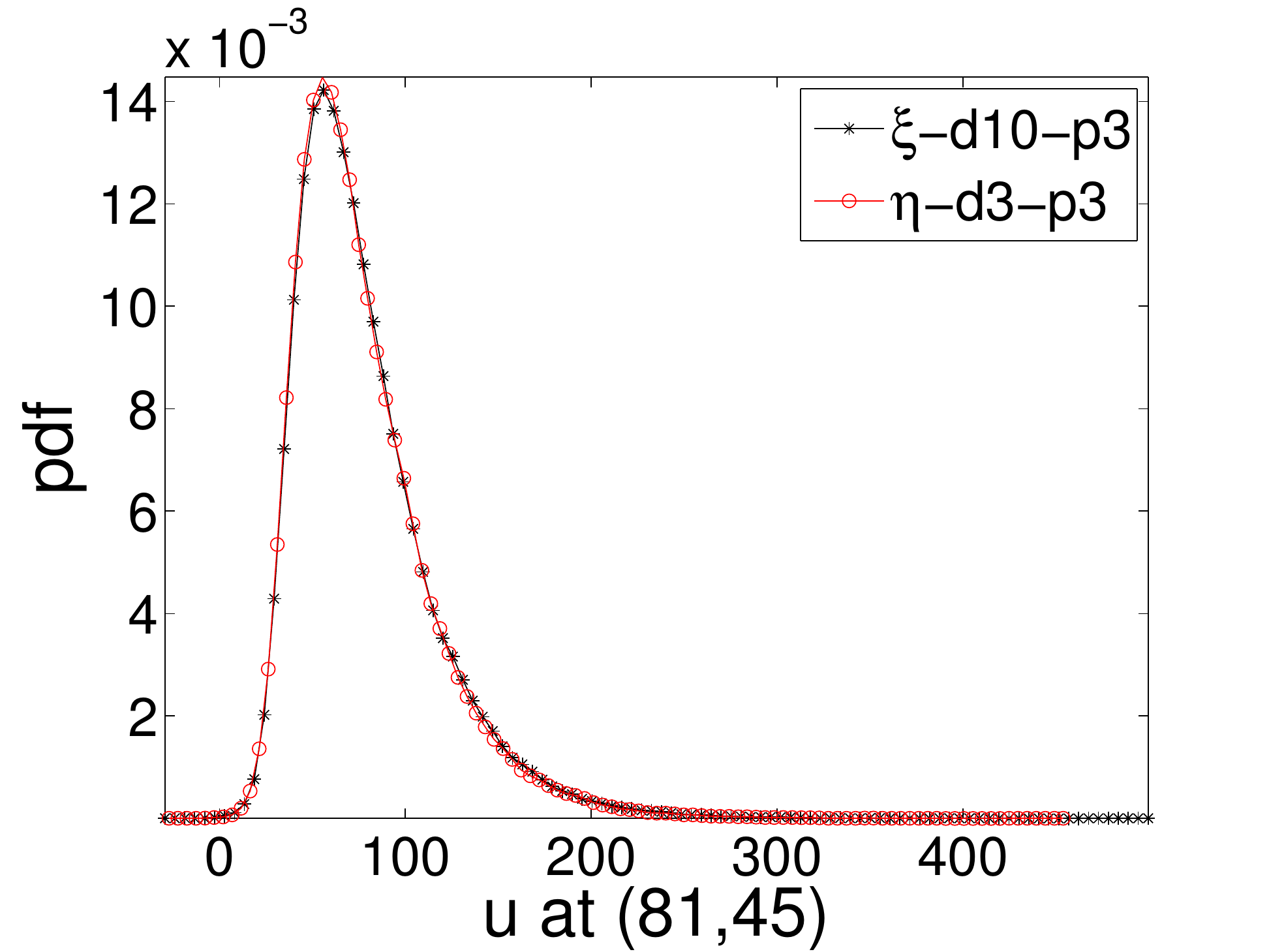}
        \caption{pdf in $D_7$, $\eta, d=3$}
     \end{subfigure}
    \begin{subfigure}[t]{0.45\textwidth}
        \centering
        \includegraphics[height=1.2in]{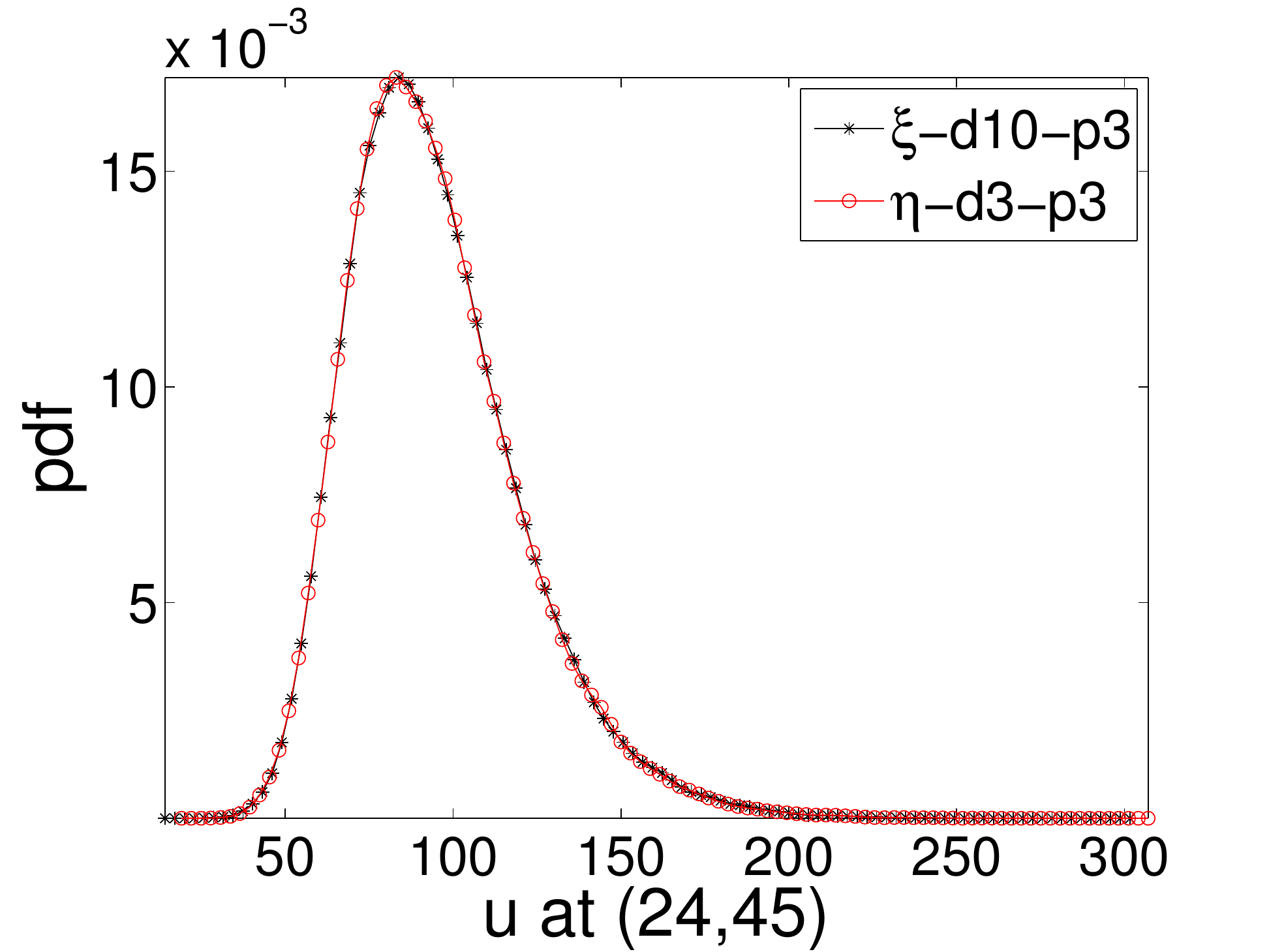}
        \caption{pdf in $D_8$, $\eta, d=3$}
     \end{subfigure}
    \caption{Probability density function of the solution at points (81,45) in $D_7$ and (24, 45) in $D_8.$} \label{RK:fig:pdf_D78}
\end{figure}

A complete characterization of the solution to a stochastic PDE is given by the probability density function (PDF), which is computed at several points in each subdomain. Figures~\ref{RK:fig:pdf_D123}, \ref{RK:fig:pdf_D456}, and \ref{RK:fig:pdf_D78} show the PDF of the solution computed in each subdomain using the full stochastic dimension and the reduced dimension adapted to that subdomain. A good agreement is observed between high- and low-dimensional PDF estimates.

\newpage

\section{Conclusions}\label{conclusions}
We have presented a novel approach to accelerate uncertainty quantification calculations for time-independent PDEs with random parameters. The main idea of the proposed approach is to decompose the spatial domain into a set of non-overlapping subdomains and use stochastic basis adaptation methods to compute a low-dimensional representation of the solution in each subdomain. We have provided numerical results where we compared the mean, standard deviation, and PDFs for solutions computed both with full and reduced dimensional representations. In our numerical experiments, the low-dimensional solutions agree quite well with high-dimensional solutions. Moreover, they do so at a significantly smaller computational cost.  

In the current work, we employ the simplest method to obtain a global solution by stitching together the solutions from the subdomains. As expected, this creates spikes in the error along subdomain boundaries. This seems from the fact that the solution in each subdomain is computed separately in a different low dimensional space. For the cases examined here, these spikes in the error are relatively small. However, they can be eliminated by using  an appropriate interpolation method or by projecting the local solution in each subdomain onto a global basis. In addition, we are working on extending the current approach to time-dependent PDEs. We will report on these issues in future publications.

%% The Appendices part is started with the command \appendix;
%% appendix sections are then done as normal sections
%% \appendix

%% \section{}
%% \label{}

%% If you have bibdatabase file and want bibtex to generate the
%% bibitems, please use

 \bibliographystyle{elsarticle-num} 
 
 \bibliography{StochDD_jcp_Rama_15Dec2015}

%% else use the following coding to input the bibitems directly in the
%% TeX file.

%\begin{thebibliography}{00}
%
%%% \bibitem{label}
%%% Text of bibliographic item
%
%\bibitem{}
%
%\end{thebibliography}

\end{document}